\documentclass[10 pt,letter,reqno]{amsart}
\usepackage{amssymb}
\usepackage{amscd}
\usepackage[all,cmtip,knot]{xy}
\xyoption{arc}

\usepackage[usenames,dvipsnames]{color}
\usepackage{hyperref}
\hypersetup{
    colorlinks=true,       			
    linkcolor=blue,          			
    citecolor=blue,		 		
    filecolor=blue,      				
    urlcolor=blue           			
}

\newtheorem*{thm}{Theorem}
\newtheorem*{prop}{Proposition}
\newtheorem*{lemma}{Lemma}

\newenvironment{pf}{\paragraph{{\sc Proof}}}{\qed\par\medskip}
\theoremstyle{definition}
\newtheorem*{defn}{Definition}
\newtheorem*{rem}{Remark}

\numberwithin{equation}{section}

\newcommand {\IC}{\mathbb{C}}
                          
\newcommand {\IP}{\mathbb{P}}
\newcommand {\IR}{\mathbb{R}}

\newcommand {\IZ}{\mathbb{Z}}

\newcommand {\A}{\mathcal A}

\newcommand {\C}{\mathcal C}
\newcommand {\D}{\mathcal D}
\newcommand {\E}{\mathcal E}

\newcommand {\J}{\mathcal J}
\newcommand {\K}{\mathcal K}

\renewcommand {\P}{\mathcal P}

\newcommand {\U}{\mathcal U}
\newcommand {\V}{\mathcal V}

\renewcommand {\b}{\mathfrak b}

\newcommand {\g}{\mathfrak{g}}
\newcommand {\h}{\mathfrak h}

\newcommand {\n}{\mathfrak n}
\newcommand {\p}{\mathfrak p}

\newcommand {\sfh}{\mathsf{h}}

\newcommand {\sfk}{\mathsf{k}}

\newcommand {\sfPhi}{\mathsf{\Phi}}

\newcommand {\ad}{\operatorname{ad}}
\newcommand {\Ad}{\operatorname{Ad}}

\newcommand {\End}{\operatorname{End}}
\newcommand {\Hom}{\operatorname{Hom}}
\renewcommand {\Im}{\operatorname{Im}}
\newcommand {\Ker}{\operatorname{Ker}}

\newcommand {\id}{\operatorname{id}}

\newcommand {\gr}{\operatorname{gr}}

\renewcommand {\sl}[1]{\mathfrak{sl}_{#1}}

\newcommand {\Ug}{U\g}

\newcommand {\Uhg}{U_{\hbar}\g}

\newcommand {\reg}{_{\operatorname{reg}}}
\newcommand {\hreg}{\h\reg}

\newcommand {\ie}{{\it i.e., }}
\newcommand {\eg}{{\it e.g.}, }
\newcommand {\fd}{finite--dimensional }

\newcommand {\lhs}{left--hand side }
\newcommand {\rhs}{right--hand side }

\newcommand {\wrt}{with respect to }

\newcommand {\ol}{\overline}
\newcommand {\wt}{\widetilde}
\newcommand {\wh}{\widehat}

\renewcommand {\DJ}{Drinfeld--Jimbo }

\newcommand {\DKKZ}{_
{\scriptscriptstyle{\operatorname{DKZ}}}}
\newcommand {\KKZ}{_
{\scriptscriptstyle{\operatorname{KZ}}}}

\newcommand {\fml}{[\negthinspace[\hbar]\negthinspace]}

\newcommand {\fmls}[1]{[\negthinspace[#1]\negthinspace]}

\newcommand {\half}[1]{\frac{#1}{2}}

\newcommand {\aand}{\qquad\text{and}\qquad}

\newcommand {\qh}{quasi--Hopf }
\newcommand {\qt}{quasitriangular }

\newcommand {\qtqha}{quasitriangular quasi--Hopf algebra }

\newcommand {\YB}{\operatorname{YB}}

\newcommand {\Alt}{\operatorname{Alt}}

\newcommand {\veps}{\varepsilon}

\newcommand {\comment}[1]{\footnote{\textcolor{blue}{#1}}}
\newcommand {\Omit}[1]{}

\newcommand {\ds}[1]{\displaystyle{#1}}

\newcommand {\Kalpha}{\mathcal K_\alpha}

\newcommand {\IH}{\mathbb H}

\newcommand {\onetwo}[1]{{#1}^{(1)}+{#1}^{(2)}}

\newcommand {\olmu}{{\ol{\mu}}}

\newcommand {\muone}{\mu^{(1)}}

\begin{document}

\renewcommand {\comment}[1]{}

\renewcommand {\DKKZ}{_{\operatorname{\scriptstyle{DKZ}}}}
\newcommand {\PL}{Poisson--Lie }
\newcommand {\PLg}{Poisson--Lie group }
\renewcommand {\YB}{Yang--Baxter }
\newcommand {\YBE}{Yang--Baxter equations }
\newcommand {\AM}{Alekseev--Meinrenken }
\newcommand {\EH}{Enriquez--Halbout }
\newcommand {\EEM}{Enriquez--Etingof--Marshall }
\newcommand {\AAM}{_{\scriptstyle{\operatorname{AM}}}}
\newcommand {\nr}{_{\scriptstyle{\operatorname{nr}}}}
\newcommand {\calS}{\mathcal S}
\newcommand {\UU}{\mathfrak U}
\renewcommand {\AA}{\mathfrak A}
\newcommand {\eps}{\epsilon}
\newcommand {\scl}[1]{\operatorname{scl}{(#1)}}
\newcommand {\sscl}{semiclassical limit }
\newcommand {\EEE}[1]{\End_{\IC\fml}(#1)}
\newcommand {\sfR}{{\mathsf R}}
\newcommand {\whh}{\wh{h}}
\newcommand {\Alg}{\mathsf{\operatorname{Alg}}}
\newcommand {\mm}{\mathfrak m}
\newcommand {\whS}{\wh{\calS}}
\newcommand {\whSm}{\wh{\calS}_\mm}

\renewcommand {\qh}{H}
\newcommand {\ch}{h}
\newcommand {\chh}{g}
\newcommand {\csol}{\gamma}
\newcommand {\qsol}{\Upsilon}
\newcommand {\halfplane}{\IH}

\newcommand {\Az}{A}
\newcommand {\cut}{c}
\newcommand {\cc}{counterclockwise }
\newcommand {\ccprod}{\stackrel{\curvearrowright}{\prod}}
\newcommand {\sect}{\sphericalangle}

\newcommand {\wtb}{\wt{\b}}
\newcommand {\wtg}{\wt{\g}}
\newcommand {\fsp}{\mathsf{fsp}}

\newcommand {\sfr}{\mathsf{r}}
\newcommand {\sff}{\mathsf{f}}
\newcommand {\sfj}{\mathsf{j}}
\newcommand{\RA}{R}
\newcommand{\Spec}{\operatorname{Spec}}

\newcommand {\nnu}{\nu}
\newcommand {\nnuc}{\nu^{\vee}}

\newcommand {\poiss}{\pi_J}
\newcommand {\adm}{\mathfrak{A}}

\newcommand {\sfml}{[\negthinspace[\sfh]\negthinspace]}

\newcommand {\LL}{\mathcal L}
\newcommand {\RR}{\mathcal R}

\newcommand {\SR}{\mathfrak{R}}

\newcommand {\cla}{_{\operatorname{cl}}}
\newcommand {\claj}{_{\operatorname{cl},J}}
\newcommand {\eclaj}{\mathsf{e}_{\operatorname{cl},J}}
\newcommand {\KKS}{Kirillov--Kostant--Souriau }

\title{Stokes phenomena, Poisson--Lie groups and quantum groups}
\author[V. Toledano Laredo]{Valerio Toledano Laredo}
\address{Department of Mathematics,
Northeastern University,
360 Huntington Avenue,
Boston MA 02115,
USA}
\email{V.ToledanoLaredo@neu.edu}
\author[X. Xu]{Xiaomeng Xu}
\address{School of Mathematical Sciences
\& Beijing International Center for Mathematical Research,
Beijing University,
No.5 Yiheyuan Road Haidian District, Beijing, P. R. China 100871}
\email{xxu@bicmr.pku.edu.cn}
\thanks{The first author was supported in part by the NSF grant
DMS--1802412. 
The second author was supported by the National Key Research and Development
Program of China (No. 2021YFA1002000) and by the National Natural Science
Foundation of China (No. 12171006).}
\date{March 2022}
\begin{abstract}
Let $\g$ be a complex semisimple Lie algebra, $\b_\pm\subset\g$ a pair
of opposite Borel subalgebras, and $\sfr\in\b_-\otimes\b_+$ the corresponding
solution of the classical \YB equations. Let $G$ be the simply--connected
\PLg corresponding to $(\g,\sfr)$, $H\subset B_\pm\subset G$ the subgroups
with Lie algebras $\h=\b_-\cap\b_+$ and $\b_\pm$, and $G^*=B_+\times_H
B_-$ the \PLg dual of $G$.
$G$--valued Stokes phenomena were used by Boalch \cite{Bo1,Bo2}
to give a canonical, analytic linearisation
of the \PLg structure on $G^*$.
$U\g$--valued Stokes phenomena were used by the first author to construct
a twist killing the KZ associator, and therefore give a transcendental
construction of the Drinfeld--Jimbo quantum group $\Uhg$ \cite{TL}.
In the present paper, we show that the former construction can be
obtained as semiclassical limit of the latter.
Along the way, we also show that the $R$--matrix of $\Uhg$ is a Stokes
matrix for the dynamical KZ equations.
\end{abstract}
\Omit{
arXiv abstract 2022 02 21
Let g be a complex semisimple Lie algebra, G the simply-connected
Poisson-Lie group corresponding to g, and G* its dual. G-valued Stokes
phenomena were used by Boalch [Bo1,Bo2] to give a canonical, analytic
linearisation of the Poisson structure on G*. Ug-valued Stokes phenomena
were used by the first author to construct a twist killing the KZ associator,
and therefore give a transcendental construction of the Drinfeld-Jimbo 
quantum group U_hg (arXiv:1601.04076). In the present paper, we show
that the former construction can be obtained as semiclassical limit of the latter.
Along the way, we also show that the R-matrix of U_hg is a Stokes
matrix for the dynamical KZ equations.
}
\maketitle

\setcounter{tocdepth}{1}
\tableofcontents
\newpage

\section{Introduction}%

\subsection{}

Let $\g$ be a complex semisimple Lie algebra and $\Uhg$ its quantised enveloping algebra.
The starting point of the present paper is the 
construction of $\Uhg$ from the
{\it dynamical} Knizhnik--Zamolodchikov (DKZ) equations obtained by the first author \cite{TL}.


Let $(\cdot,\cdot)$ be an invariant inner product on $\g$, $\Omega\in\g\otimes\g$ the corresponding
Casimir element, and $\h\subset\g$ a Cartan subalgebra. Consider the DKZ on $n=2$ points, that
is the $\End(\Ug^{\otimes 2})$--valued connection on $\IC\ni z=z_1-z_2$ given by
\begin{equation}\label{eq:intro DKZ}
d-\left(\sfh\frac{\Omega}{z}+\ad\muone\right) dz
\end{equation}
where $\mu\in\h$, $\muone=\mu\otimes 1$, and $\sfh$ is a formal deformation parameter. Just as
its non--dynamical counterpart which is obtained for $\mu=0$, the connection \eqref{eq:intro DKZ}
has a regular singularity at $z=0$, and admits a canonical fundamental solution $\Upsilon_0$ which
is asymptotic to $z^{\sfh\Omega}$ as $z\to 0$.

\subsection{}\label{eq:intro J}

The dynamical term $\ad\muone$ gives rise to an {\it irregular singularity} at $z=\infty$. Assuming that
$\mu$ is real, so that all Stokes rays lie in $\IR$, and regular, it is proved in \cite{TL} that \eqref{eq:intro DKZ}
admits two canonical fundamental solutions $\Upsilon_\pm$ which are asymptotic to $e^{z\ad\muone}
\cdot z^{\sfh\Omega_0}$ as $z\to\infty$ with $\Im z\gtrless 0$, where $\Omega_0\in\h\otimes\h$ is the
projection of $\Omega$.

Consider now the regularised holonomy of \eqref{eq:intro DKZ} from $\pm\iota\infty$ to $0$ \ie the
element $J_\pm\in\Ug^{\otimes 2}\sfml$ given by $J_\pm=\Upsilon_0(z)^{-1}\cdot\Upsilon_\pm(z)$,
where $\Im z\gtrless 0$. One of the main results of \cite{TL} is that $J_\pm$, regarded as a twist,
kills the KZ associator $\Phi\KKZ$ which arises from the (non--dynamical, reduced) KZ equations
on $n=3$ points
\[d-\sfh\left(\frac{\Omega_{12}}{z}+\frac{\Omega_{23}}{z-1}\right)dz\]

Let $\Delta_\pm=J_\pm^{-1}\Delta(\cdot)J_\pm$ and
$R_\pm=(J_\pm^{21})^{-1}e^{\hbar\Omega/2}J_\pm$ be the corresponding twisted coproduct and
$R$--matrix, where $\hbar=2\pi\iota\sfh$. It follows that $\left(\Ug\fml,\Delta_\pm,R_\pm\right)$ is
a quasitriangular Hopf algebra, which can be shown to be isomorphic to the quantum group $\Uhg$.

\subsection{}

In contrast to earlier constructions of $\Uhg$ from the (non--dynamical) KZ equations
\cite{Dr4,KL,ek-1}, the above construction is entirely transcendental \ie does not rely
on cohomological arguments or the representation theory of $\g$, and perhaps more
naturally explains how $\Uhg$ arises from such equations.

One additional feature is its compatibility with the Casimir equations of $\g$ introduced in \cite{DC,
MTL,TL00,FMTV}. Specifically, the twist $J_\pm$ is a smooth function of $\mu\in\hreg^{\IR}$, and
satisfies the PDE
\[d J_\pm=
\half{\sfh}\sum_{\alpha\in\sfPhi_+}\frac{d\alpha}{\alpha}
\left(\Delta(\Kalpha)J_\pm-J_\pm(\onetwo{\Kalpha})\right)\]
where $\sfPhi_+$ is a chosen system of positive roots, and $\Kalpha$ the Casimir of the $\sl{2}
$--subalgebra of $\g$ corresponding to $\alpha$. This compatibility is a key ingredient in proving
that the monodromy of the Casimir connection of $\g$ is given by Lustzig's quantum Weyl group
operators \cite{TL1,TL,ATL3}. 


\subsection{} 

Let now $G$ be the connected and simply connected complex Lie group corresponding to $\g$.
Irregular singularities were exploited earlier by Boalch to linearise the Poisson structure on the
\PLg $G^*$ dual to $G$ \cite{Bo1,Bo2}.

Boalch considered connections on the holomorphically trivial $G$--bundle over $\IP^1$ of the form
\begin{equation}\label{eq:intro IODE}
d-\left(\frac{A}{z^2}+\frac{B}{z}\right) dz
\end{equation}
where $A\in\h$ is regular, and $B\in\g$.
%
%
Assume that $A$ is real, so that the Stokes rays of \eqref{eq:intro IODE} lie in $\IR$,\footnote
{Contrary to the reality assumption made in \ref{eq:intro J}, the assumption that $A\in\h^\IR$
is inessential, and is only made in the Introduction to simplify the exposition.} and set $\IH_
\pm=\{z\in\IC|\,\Im z\gtrless 0\}$. Then, there are unique holomorphic fundamental solutions
$\gamma_\pm:\IH_\pm\to G$ of \eqref{eq:intro IODE}, which are asymptotic to $e^{-A/z}\cdot
z^{[B]}$ as $z\to 0$ in $\IH_\pm$, where $[B]$ is the projection of $B$ onto $\h$.

Define the {\it Stokes matrices} $S_\pm\in G$ by the analytic continuation identities
\[ \wt{\gamma_-}=\gamma_+\cdot S_+
\aand
\wt{\gamma_+}=
\gamma_-\cdot S_-\cdot e^{2\pi\iota[B]}\]
where $\wt{\cdot}$ denotes counterclockwise analytic continuation, and the identities hold
in $\IH_+$ and $\IH_-$ respectively. The elements $S_\pm\in G$ are unipotent. Specifically,
$A$ determines a partition $\sfPhi=\sfPhi_+\sqcup\sfPhi_-$ of the root system by $\sfPhi_
\pm=\{\alpha\in\sfPhi|\,\alpha(A)\gtrless 0\}$, and $S_\pm$ lies in the unipotent subgroup
$N_\pm\subset G$ with Lie algebra $\n_\pm=\bigoplus_{\alpha\in\sfPhi_\pm}\g_\alpha$.

\subsection{} 

Let $B_\pm\subset G$ be the Borel subgroups corresponding to $\sfPhi_\pm$,
$H=B_+\cap B_-$ the maximal torus with Lie algebra $\h$, and
consider the fibred product
\[B_+\times_H B_-=\{(b_+,b_-)\in B_+\times B_-|\, \pi_+(b_+)\pi_-(b_-)=1\}\]
where $\pi_\pm:B_\pm\to H$ are the quotient maps. Following \cite
{Bo2}, we define the {\it Stokes map} to be the analytic map $\calS:
\g\longrightarrow B_+\times_H B_-$ given by
\begin{equation}\label{eq:intro Stokes}
B\longrightarrow
\left(S_+^{-1}\cdot e^{-\iota\pi[B]},
S_-\cdot e^{\iota\pi[B]}\right)
\end{equation}

\Omit{Note that $B_-\times_H B_+$ maps to $G$ via the map 
$\beta:(b_+,b_-)\to b_+\cdot b_-^{-1}$.\comment{$\beta$ is a principal
bundle over its image with structure group the order two elements in
$H$.} Moreover, by Proposition \ref{pr:monodromy reln}, the composition
$\beta\circ\calS_r$ is the map $\g\to G$ given by 
\[B\longrightarrow 
(S^r_+)^{-1}\cdot e^{-2\pi\iota[B]}\cdot (S^r_-)^{-1}=
C_r\cdot e^{-2\pi\iota B}\cdot C_r^{-1}\]}

\subsection{} 

The pair $(B_+,B_-)$ 
gives rise to a solution $\sfr\in\b_-\otimes
\b_+$ of the classical \YBE given by
\begin{equation}\label{eq:intro standard cybe}
\sfr = x_i\otimes x^i + \half{1} t_a\otimes t^a
\end{equation}
where $\{x_i\},\{x^i\}$ are bases of $\n_-,\n_+$ which are
dual \wrt $(\cdot,\cdot)$, and $\{t_a\},\{t^a\}$ are dual bases of $\h$.
The element $\sfr$ gives $\g$ the structure of a \qt Lie bialgebra, with
cobracket $\delta:\g\to\g\wedge\g$ given by $\delta(x)=[x\otimes 1+1
\otimes x,\sfr]$. 

The dual Lie bialgebra $(\g^*,\delta^t,[\cdot,\cdot]^t)$ may be identified, as
a Lie algebra, with
\[\b_+\times_\h\b_-=\{(X_+,X_-)\in\b_+\oplus\b_-|\pi_+(X_+)+\pi_-(X_-)=0\}\]
where $\pi_\pm:\b_\pm\to\h$ is the quotient map.
This endows $G^*=B_+\times_H B_-$ with the structure of a \PL group,
which is dual to $G$.

\subsection{} 

Endow $\g^*$ with its standard \KKS Poisson structure 
\[ \{f,g\}(x)=\langle[d_x f,d_xg],x\rangle\]
where $d_xh\in T^*_x\g^*=\g$ is the differential of $h$ at $x$, 
and $[\cdot, \cdot]$ is the Lie bracket on $\g$.

Let $\nnu:\g^*\to\g$ be the isomorphism induced by the bilinear form
$(\cdot,\cdot)$, and identify $\g$ and $\g^*$ by using $\nnu^\vee=-1/
(2\pi\iota)\nnu$. The following remarkable result is due to Boalch \cite
{Bo1,Bo2}.

\begin{thm}
The map $\calS:\g^*\to G^*$ is a Poisson map, and generically a local
complex analytic diffeomorphism. In particular, $\calS$ gives a linearisation
of the Poisson structure on $G^*$. 
\end{thm}

\Omit{
Let $\nnu:\g^*\to\g$ be the isomorphism induced by the bilinear form
$(\cdot,\cdot)$, and set $\nnu^\vee=-1/(2\pi\iota)\nnu$. The following
remarkable result is due to Boalch \cite{Bo1,Bo2}.

\begin{thm}\label{th:Boalch}
The map $\calS\circ\nnu^\vee:\g^*\to G^*$ is a Poisson map, and 
generically a local complex analytic diffeomorphism.
In particular, $\calS\circ\nnu^\vee$ gives a linearisation of the
Poisson structure on $G^*$. 
\end{thm}
}



\subsection{}

One of goals of the present paper is to prove that Boalch's linearisation result, specifically
the fact that $\calS$ is a Poisson map, can be obtained as a semiclassical limit of the
transcendental construction of $\Uhg$.

Our overall strategy is the following. Since $\calS$ is holomorphic, it suffices to show that its
Taylor series $\whS$ at $0\in\g^*$ is a formal Poisson map. This in turn follows if $\whS$ can
be quantised. We therefore seek quantisations $\IC_\hbar\fmls{\g^*}$ and $\IC_\hbar\fmls{G^*}$
of the algebras of functions on the formal \PL groups corresponding to $\g^*$ and $G^*$, together
with an algebra isomorphism $\whS^*_\hbar:\IC_\hbar\fmls{G^*}\to\IC_\hbar\fmls{\g^*}$ such
that the following diagram is commutative
\begin{equation}\label{eq:diagram}
\xymatrix@C=2cm{
\IC_\hbar\fmls{\g^*}\ar[d] 	& \IC_\hbar\fmls{G^*}\ar[l]_{\whS^*_\hbar}\ar[d]\\
\IC\fmls{\g^*}				& \IC\fmls{G^*}\ar[l]^{\whS^*}\\
}\end{equation}
where the vertical arrows are the specialisations at $\hbar=0$, and the bottom one the 
pullback of $\whS$.

\subsection{}\label{ss:intro duality}

A formal quantisation of the dual $P^*$ of a \PLg $P$ can be obtained from Drinfeld's {\it quantum
duality principle} as follows \cite{DrICM,Ga}. Let $\UU$ be a quantised enveloping algebra which
deforms the Lie bialgebra $\p$ of $P$. Thus, $\UU$  is a topologically free Hopf algebra over $
\IC\fml$ such that $\UU/\hbar\UU$ is isomorphic to $U\p$ and, for any $x\in\p$ with cobracket
$\delta(x)\in\p\wedge\p$
\[\delta(x)=\left.\frac{\Delta(\wt{x})-\Delta^{21}(\wt x)}{\hbar}\right|_{\hbar=0}\]
where $\wt{x}\in\UU$ is an arbitrary lift of $x$. Then, $\UU$ admits a canonical Hopf subalgebra
$\UU'$ which is commutative mod $\hbar$, and endowed with a canonical Poisson isomorphism
$\imath_\UU:\UU'/\hbar\UU'\to\IC\fmls{P^*}$.

The simplest example of Drinfeld duality arises when $P$ is the Lie group $G$ endowed with
the trivial Poisson structure. The corresponding Lie bialgebra is $\g$ with the trivial cobracket,
and $P^*$ is the additive abelian group $\g^*$ with cobracket given by the transpose of the bracket
on $\g$. In this case, $\UU$ can be taken to be $\Ug\fml$ with undeformed product and coproduct.
The corresponding subalgebra $\UU'$ is the Rees algebra of formal power series $\sum_{n\geq 0}
x_n\hbar^n$ where the filtration order of $x_n$ is at most $n$, and the specialisation $\imath_
\UU$ is the symbol map $\UU'/\hbar\UU'\to\prod_{n\geq 0}S^n\g=\IC\fmls{\g^*}$. 

\subsection{}


To obtain a formal quantisation of $G^*$, we seek a QUE deforming the quasitriangular
Lie bialgebra $(\g,\sfr)$, where $\sfr\in\b_-\otimes\b_+$ is the canonical element \eqref
{eq:intro standard cybe}. 
One such quantisation is the
\DJ quantum group $\Uhg$ corresponding to $\g$. That, however, shifts the problem of
filling in the diagram \eqref{eq:diagram} to one of finding an algebra isomorphism $(\Uhg)'
\to\UU'$, where $\UU=\Ug\fml$, and showing that the latter quantises $\whS^*$.

Alternatively, we may resort to a {\it preferred} quantisation of $\g$, that is a QUE
which is equal to $\UU$ as algebras. A class of such quantisations may be obtained
as a {\it twist quantisation}, that is by using an element $J\in 1+\frac{\hbar}{2}\sfj+
\hbar^2\UU^{\otimes 2}$ satisfying $\sfj-\sfj^{21}=\sfr-\sfr^{21}$, together with the
twist equation
\[\Phi\cdot J_{12,3}\cdot J_{1,2}=J_{1,23}\cdot J_{2,3}\]
where $\Phi$ is a given associator. Then, 
$\UU_J=
\left(\UU,J^{-1}\Delta_0(\cdot)J,J^{-1}_{21}e^{\hbar\Omega/2}J\right)$
is a QUE which quantises $(\g,\sfr)$, and $(\UU_J)'$ is a formal quantisation
of $G^*$. 

\subsection{} 

A general result of \EH asserts that if the twist $J$ is {\it admissible}, that is such
that $\hbar\log(J)\in(\UU')^{\otimes 2}$, the Drinfeld algebras $(\UU_J)'$ and $\UU'$
{\it coincide} \cite{EH}. In this case, the equality $e:(\UU_J)'\to\UU'$ clearly is an
algebra isomorphism, and descends to a Poisson isomorphism
$e\claj:\IC\fmls{G^*}\to\IC\fmls{\g^*}$ given by the composition
\[e\claj = \imath_\UU \circ e_0 \circ\imath_{\UU_J}^{-1}\]
where $\imath_\UU:\UU'/\hbar\UU'\to\IC\fmls{\g^*}$ is the symbol map, $\imath_
{\UU_J}:(\UU_J)'/\hbar(\UU_J)'\to\IC\fmls{G^*}$ the canonical identification mentioned
in \ref{ss:intro duality}, and $e_0=\id$ the reduction of $e$ mod $\hbar$.

One of the main results of this paper is that if $J=J_+$ is (one of) the twist(s)
arising from the dynamical KZ equations described in \ref{eq:intro J}, 
then $J$ is admissible, and the corresponding map $e\claj$ is equal to the
Stokes map $\whS^*$. In particular, the latter is a Poisson map.

\subsection{} \label{ss:EEM}

A key ingredient in proving the identity $e\claj=\whS^*$ is a result of \EEM \cite{EEM}
which gives an explicit formula for $e\claj$, under the additional assumptions that
$\Phi$ is a Lie associator and that the admissible twist $J$ lies in $\UU'\otimes\UU
\cap\UU\otimes\UU'$.

Consider to that end the quotient $\UU\otimes\UU'/\hbar\,\UU\otimes\UU'\cong\Ug
\fmls{\g^*}$, where the latter is the algebra of $\Ug$--valued formal power series
on $\g^*$. Let $G\fmls{\g^*}_+\subset\Ug\fmls{\g^*}$ be the prounipotent group of
$\IC\fmls{\g^*}$--points of 
$G$ such that their value
at $0\in\g^*$ is equal to $1$. Then, the following holds \cite{EEM}
\begin{enumerate}
\item The {\it semiclassical limit $\jmath=\scl{J}$}, that is the image of $J$ in $\Ug\fmls{\g^*}$,
lies in $G\fmls{\g^*}_+$ and is therefore a formal map $\g^*\to G$.\\
\item 
Let
\[\beta:G^*\to G,\qquad (b_+,b_-)\to b_+\cdot b_-^{-1}\]
be the {\it big cell map}. Then, the composition of $e\claj$ with $\beta$
is the formal map $\g^*\to G$ given by the {\it twisted exponential map}
\begin{equation}\label{eq:intro eem}
e_\jmath(\lambda)=
\jmath(\lambda)^{-1}\cdot e^{\nnu(\lambda)}\cdot\jmath(\lambda)
\end{equation}
where $\nu:\g^*\to\g$ is the isomorphism given by the inner product.
\end{enumerate}

\subsection{} 

Since $\beta$ is an isomorphism when regarded as a formal map, it suffices to prove
that $e\claj\circ\wh{\beta}^*=\whS^*\circ\wh{\beta}^*=\wh{\beta\circ\calS}^*$ that is, by
\eqref{eq:intro eem} that
\begin{equation}\label{eq:intro to prove}
\wh{\beta\circ\calS}=\jmath(\lambda)^{-1}\cdot e^{\nnu(\lambda)}\cdot\jmath(\lambda)
\end{equation}

By definition of $\calS$ \eqref{eq:intro Stokes}, the composition $\beta\circ\calS$
is the map $B\to S_+^{-1}\cdot e^{-2\pi\iota[B]}\cdot S_-^{-1}$, which is the {\it 
clockwise} monodromy around $z=0$ of \eqref{eq:intro IODE} expressed in the 
solution $\gamma_-$. By parallel transport to $z=\infty$, where \eqref{eq:intro IODE}
has a regular singularity with residue $-B$, $\beta\circ\calS$ is also equal to
\begin{equation}\label{eq:intro infty}
B\to C^{-1}_-\cdot e^{-2\pi\iota B}\cdot C_-
\end{equation}
where $C_-=C_-(B)\in G$ is the {\it connection matrix}, 
that is the element relating $\gamma_-$ to the canonical
fundamental solution $\gamma_\infty$ which is asymptotic to $z^
{2\pi\iota B}$ near $z=\infty$.
\Omit{
where $C_-=C_-(B)\in G$ is the {\it connection matrix} of \eqref{eq:intro IODE} 
that is the element defined by $\gamma_-(z)=\gamma_\infty(z)\cdot C_-$, $z\in
\IH_-$, where $\gamma_\infty$ is the canonical fundamental solution asymptotic
to $z^{2\pi\iota B}$ near $z=\infty$.
}

Comparing the right--hand sides of \eqref{eq:intro to prove} and \eqref{eq:intro infty},
and recalling that $\g^*$ and $\g$ are identified by $-1/2\pi\iota\cdot\nu$, it therefore
suffices to show that $B\to\wh{C_-}$ is the semiclassical limit of the DKZ twist $J$.\footnote
{The problem of obtaining a quantisation of the connection matrix $C_-$ formulated
in \cite{Xu}, together with our intuition that such a quantisation should be given by
the DKZ twist $J$, were in fact the original impetus of this project.}

\subsection{}

The fact that $J$ is a quantisation of the connection matrix $C_-$ follows from
the uniqueness of canonical fundamental solutions of 
\eqref{eq:intro IODE}, when the structure group is an arbitrary affine algebraic group,
specifically the prounipotent group $G\fmls{\g^*}_+$ \cite{BTL2}.
It stems from the basic, but seemingly novel observation that the {\it semiclassical
limit of the DKZ equation \eqref{eq:intro DKZ} is equal to the ODE \eqref
{eq:intro IODE}}.\footnote{This is related to, but different from, the fact that
a different semiclassical limit of the KZ equations are the (non--linear)
Schlesinger equations \cite{Re}.}

More precisely, if $\Upsilon$ is a solution of
\[\frac{d\Upsilon}{dz}=
\left(\ad\muone+\sfh\frac{\Omega}{z}\right)\Upsilon\]
with values in $\UU\otimes\UU'$, the semiclassical limit $\gamma$
of $\Upsilon$, as a formal function of $\lambda\in\g^*$ with values
in $\Ug$, is readily seen to satisfy
\[\frac{d\gamma}{dz}=
\left(\ad\mu+\frac{\nnu(\lambda)}{2\pi\iota z}\right)\gamma\]
where $\nnu(\lambda)=\id\otimes\lambda(\Omega)$ which, after
the change of variable $z\to 1/z$, and the replacement $\ad\mu
\to -\Az, \nnu(\lambda)\to -2\pi\iota B$ is precisely the equation \eqref
{eq:intro IODE}.\footnote{The appearance of the factor $2\pi\iota$
is due to the fact that the identification $\UU'/\hbar\UU'\cong\wh{S\g}$
is given by mapping $x\in\g$ to $\hbar x=2\pi\iota\sfh x\in\UU'$.}

\subsection{Outline of paper} 

In Sections \ref{se:G Stokes} and \ref{se:Ug Stokes}, we review the definition of
the Stokes data and map for the connection \eqref{eq:intro IODE}, and the 
transcendental construction of $\Uhg$ given in \cite{TL}. In Section \ref{se:R = S}, 
we show that quantum $R$--matrix of $\Uhg$ is a Stokes matrix of the dynamical
KZ equation. Section \ref{se:quantum duality} reviews Drinfeld's duality principle.
Section \ref{se:scl DKZ} contains the first part of our main results, namely the fact that
the semiclassical limits of the DKZ equations and its canonical solutions at $0$ and
$\infty$ are equal to the connection \eqref{eq:intro IODE} and its canonical solutions,
after a change of variables. 
Section \ref{se:EEM} describes the linearisation formula of Enriquez--Etingof--Marshall.
Finally, in Section \ref{se:alt boalch}, we prove that the Stokes map is Poisson and,
in Section \ref{Section:isomono} relate quantum and classical isomonodromic equations.

\subsection*{Acknowledgements}
\noindent
We would like to thank Anton Alekseev and Pavel Etingof for their helpful discussions and useful comments.
\section{Stokes phenomena and \PL groups} \label{se:G Stokes}

\subsection{$G$--valued irregular connections on $\IP^1$}\label{ss:G Stokes}


Let $G$ be an affine algebraic group defined over $\IC$, $H\subset G$
a maximal torus, and $\h\subset\g$ the Lie algebras of $H$ and $G$
respectively. Let $\sfPhi\subset\h^*$ be the set of roots of $\g$ relative
to $\h$, and $\hreg=\h\setminus\bigcup_{\alpha\in\sfPhi}\Ker\alpha$
the set of regular elements in $\h$. 

Let $\P$ be the holomorphically trivial, principal $G$--bundle on
$\IP^1$, and consider the meromorphic connection $\nabla$ on
$\P$ given by
\begin{equation}\label{eq:dKZ}
\nabla=d-\left(\frac{\Az}{z^2}+\frac{B}{z}\right)dz.
\end{equation}
where $\Az,B\in\g$. We assume henceforth that $\Az\in\hreg$.
By definition, the {\it Stokes rays} of $\nabla$ are the rays $\IR
_{>0}\cdot\alpha(\Az)\subset\IC^*$, $\alpha\in\sfPhi$, that is the
rays through the non--zero eigenvalues of $\ad(A)$. 
A ray $r$ is called {\it admissible} if it is not a Stokes ray.

\subsection{Canonical fundamental solutions}

To each admissible ray $r$, and determination of $\log z$, there
is a canonical fundamental solution $\csol_r$ of $\nabla$ with
prescribed asymptotics in the open half--plane
\[\halfplane_r=\left\{ue^{\iota\phi}|\,u\in r, \phi\in(-\pi/2,\pi/2)\right\}\]
Specifically, the following result is proved in \cite{BJL} for $G=GL
_n(\IC)$, in \cite{Bo2} for $G$ reductive, and in \cite{BTL2} for an
arbitrary affine algebraic group.\footnote{We use the formulation
of \cite{BTL2}, which does not rely on formal power series solutions.}
Denote by $[B]$ the projection of $B$ onto $\h$ corresponding to
the root space decomposition $\g=\h\bigoplus_{\alpha\in\sfPhi}\g_\alpha$.

\begin{thm}\label{jurk}
Let $r=\IR_{>0}\cdot e^{\iota\theta}$ be an admissible ray. Then,
there is a unique holomorphic function $\ch_r:\halfplane_r\to G$
such that
\begin{enumerate}
\item $\ch_r$ tends to 1 as $z\to 0$ in any closed sector of
$\halfplane_r$ of the form
\[|\arg(z\cdot e^{-\iota\theta})|\leq \frac{\pi}{2}-\delta,\qquad\delta>0\]
\item For any determination of $\log z$ with a cut along the ray $c$,
the function
\[\csol_r=\ch_r\cdot e^{-\Az/z}\cdot z^{[B]}\]
where $z^{[B]}=\exp([B]\log z)$, satisfies $\nabla\csol_r=0$ on
$\IH_r\setminus c$.
\end{enumerate}
\end{thm}

\subsection{Stokes phenomena}

For a given determination of $\log z$, with a cut along a ray $\cut$,
the canonical solution $\csol_r$ is locally constant \wrt the choice
of $r$, so long as $r$ does not cross a Stokes ray. More precisely,
the following holds. For any subset $\Sigma\subset\IC$, let $\g_
\Sigma\subseteq\g$ be the direct sum of the eigenspaces of
$\ad(A)$ corresponding to the eigenvalues contained in $\Sigma$,
\[\g_\Sigma=\bigoplus_{\substack{\alpha\in\sfPhi\sqcup\{0\}:\\ \alpha(\Az)\in\Sigma}}\g_\alpha\]
where $\g_0=\h$. Note that $[\g_{\Sigma_1},\g_{\Sigma_2}]
\subseteq\g_{\Sigma_1+\Sigma_2}$. In particular, if $\Sigma$ is
an open convex cone, $\g_\Sigma$ is a nilpotent subalgebra of
$\g$.

\begin{prop}\label{pr:r and r'}
Let $r,r'$ be admissible rays such that $r\neq-r'$, so that $\halfplane
_r\cap\halfplane_{r'}\neq\emptyset$, and denote by $\ol{\Sigma}(r,r')
\subset\IC^\times$ the closed convex cone bounded by $r$ and $r'$.
Let
\[S:\halfplane_r\cap\halfplane_{r'}\setminus\cut\longrightarrow G\]
be the locally constant function defined by $\csol_r=\csol_{r'}\cdot S$.
Then, the following holds.
\begin{enumerate}
\item $S$ takes values in the unipotent elements of $G$,
and $\log S$ in the nilpotent subalgebra $\g_{\ol{\Sigma}
(r,r')}$.
\item In particular, if $\ol{\Sigma}(r,r')$
does not contain any Stokes rays, the solutions $\csol_r$ and
$\csol_{r'}$ coincide on $\halfplane_r\cap\halfplane_{r'}\setminus\cut$.
\end{enumerate}
\end{prop}
\begin{pf}
The asymptotic behaviour of $\csol_r$ and $\csol_{r'}$ implies that
\begin{equation}\label{eq:S asy}
z^{[B]}\cdot e^{-\Az/z}\cdot S\cdot e^{\Az/z}\cdot z^{-[B]}
=
\left(\gamma_{r'}\cdot e^{\Az/z}\cdot z^{-[B]}\right)^{-1}\cdot\gamma_r\cdot e^{\Az/z}\cdot z^{-[B]}
\to 1
\end{equation}
as $z\to 0$ along any ray $\rho$ in $\halfplane_r\cap \halfplane
_{r'}\setminus\cut$.
By \cite[Lemma 6]{Bo2} and \cite[Prop. 6.3]{BTL2}, the restriction
of $S$ to $\rho$ is unipotent, and $\log S$ lies in $\g_{\halfplane
_\rho}$.

Up to a permutation, we may assume that the counterclockwise
angle from $r$ to $r'$ is less than $\pi$, so that $\IH_r\cap\IH_{r'}$
is the open convex cone bound by the rays $r'e^{-i\pi/2}$ and $re^
{i\pi/2}$.

If the cut $\cut$ is not contained in $\halfplane_r\cap\halfplane_{r'}$,
$S$ takes a single value. Since the intersection of the half--planes
$\halfplane_\rho$ as $\rho$ varies in $\halfplane_r\cap\halfplane_{r'}$
is the closed convex cone bounded by $r$ and $r'$, it follows that
$\log S\in\g_{\ol{\Sigma}(r,r')}$. 

If, on the other hand, $\cut$ disconnects $\halfplane_r\cap\halfplane
_{r'}$ into two open cones $\Sigma_<,\Sigma_>$, listed in counterclockwise
order, then $\csol_r=\csol_{r'}\cdot S_\lessgtr$ on $\Sigma_\lessgtr$,
for some $S_\lessgtr\in G$. 
The previous argument
then shows that $S_\lessgtr$ are unipotent, and that
\[\log S_<\in\g_{\ol{\Sigma}(e^{-\iota\pi/2}c,r')}
\aand
\log S_>\in\g_{\ol{\Sigma}(r,e^{\iota\pi/2}c)}\]
Analytic continuation across $\cut$ implies that $S_>=e^{2\pi\iota[B]}
\cdot S_<\cdot e^{-2\pi\iota[B]}$. Since any $\g_\Sigma$ is stable under
$\Ad(e^{2\pi\iota[B]})$, this implies that
\[\log S_\lessgtr\in
\g_{\ol{\Sigma}(e^{-\iota\pi/2}c,r')}\cap
\g_{\ol{\Sigma}(r,e^{\iota\pi/2}c)}
=\g_{\ol{\Sigma}(r,r')}\]
\end{pf}

\subsection{Stokes data}

For any two rays $r,r'$, let $\sect(r,r')\subset\IC^\times$ be the (not
necessarily convex) closed sector swept by $e^{\iota\theta}\cdot r$,
as $\theta$ ranges from $0$ to the positive angle between $r$ and
$r'$. If $r,r'$ are admissible, and different from the log cut $\cut$,
define an element $S_{r'r}\in G$ by the identity
\[\left.\wt{\left.\csol_{r}\right|_{r}}\right|_{r'}=
\left.\csol_{r'}\right|_{r'}\cdot S_{r',r}\cdot e^{2\pi\iota[B]\epsilon^\cut_{r'r}}\]
where the \lhs is the \cc analytic continuation to $r'$ of the restriction
of $\csol_r$ to $r$, and $\epsilon_{r'r}^\cut$ is $1$ if $\cut$ lies
in $\sect(r,r')$, and $0$ otherwise. 

\begin{prop}\label{pr:uni and fact}
The following holds
\begin{enumerate}
\item If the positive angle formed by $r$ and $r'$ is at most $\pi$, $S_{r'r}$
is unipotent, and its logarithm lies in the nilpotent subalgebra $\g_{\sect(r,r')}$.
\item If the admissible ray $r'\neq c$ lies in $\sect(r,r'')$, the following
factorisation holds
\[ S_{r''r}=
S_{r''r'}\cdot e^{2\pi\iota[B]\epsilon^\cut_{r''r'}}\cdot S_{r'r}\cdot e^{-2\pi\iota[B]\epsilon^\cut_{r''r'}}\]
\end{enumerate}
\end{prop}
\begin{pf}
(1) Let $\ell$ be a ray in $\IH_r\cap\IH_{r'}\setminus\cut$. Then
\[\left.\wt{\left.\csol_r\right|_r}\right|_\ell=\left.\csol_r\right|_\ell\cdot e^{2\pi\iota[B]\cdot\epsilon^\cut_{\ell r}}
\aand
\left.\wt{\left.\csol_{r'}\right|_\ell}\right|_{r'}=\left.\csol_{r'}\right|_{r'}\cdot e^{2\pi\iota[B]\cdot\epsilon^\cut_{r'\ell}}\]
By Proposition \ref{pr:r and r'}, $\left.\csol_r\right|_{\ell}=\left.\csol_{r'}
\right|_{\ell}\cdot S$, where $S\in G$ is a unipotent element whose log
lies in $\g_{\ol{\Sigma}(r,r')}$. Computing analytic continuation in stages
yields
\[\begin{split}
\left.\wt{\left.\csol_r\right|_r}\right|_{r'}
&=
\left.\wt{\left.\wt{\left.\csol_r\right|_r}\right|_{\ell}}\right|_{r'}
= 
\left.\wt{\left.\csol_r\right|_\ell}\right|_{r'}\cdot e^{2\pi\iota[B]\cdot\epsilon^\cut_{\ell r}}\\
&= 
\left.\wt{\left.\csol_{r'}\right|_\ell}\right|_{r'}\cdot S\cdot e^{2\pi\iota[B]\cdot\epsilon^\cut_{\ell r}}
= 
\left.\csol_{r'}\right|_{r'}\cdot e^{2\pi\iota[B]\cdot\epsilon^\cut_{r'\ell}}\cdot S\cdot e^{2\pi\iota[B]\cdot\epsilon^\cut_{\ell r}}\\
&= 
\left.\csol_{r'}\right|_{r'}\cdot \Ad(e^{2\pi\iota[B]\cdot\epsilon^\cut_{r'\ell}})(S)\cdot e^{2\pi\iota[B]\cdot\epsilon^\cut_{r'r}}
\end{split}\]
so that $S_{r'r}=\Ad(e^{2\pi\iota[B]\cdot\epsilon^\cut_{r'\ell}})(S)$.

(2) Computing analytic continuation from $r$ to $r''$ in stages yields
\[\begin{split}
\left.\wt{\left.\csol_r\right|_r}\right|_{r''}
&= 
\left.\wt{\left.\wt{\left.\csol_r\right|_r}\right|_{r'}}\right|_{r''}
=
\left.\wt{\left.\csol_{r'}\right|_{r'}}\right|_{r''}\cdot S_{r'r}\cdot e^{2\pi\iota[B]\epsilon^\cut_{r'r}}\\
&=
\left.\csol_{r''}\right|_{r''}\cdot S_{r''r'}\cdot e^{2\pi\iota[B]\epsilon^\cut_{r''r'}}\cdot S_{r'r}\cdot e^{2\pi\iota[B]\epsilon^\cut_{r'r}}
\end{split}\]
Since the result is also equal to $\left.\csol_{r''}\right|_{r''}\cdot S_{r''r}\cdot
e^{2\pi\iota[B]\epsilon^\cut_{r''r}}$, the result follows.
\end{pf}

\subsection{Stokes factors}

Given a Stokes ray $\ell$, the {\it Stokes factor} $S_\ell$ is the
unipotent element of $G$ defined by $S_\ell=S_{r'r}$, where $
r,r'\neq\cut$ are admissible rays such that $\sect(r,r')$ contains
no other Stokes rays than $\ell$, and does not contain the cut
$\cut$ if the latter is different from $\ell$. By Proposition \ref
{pr:uni and fact}, the definition of $S_\ell$ is independent of
the choice of $r,r'$. The following is a direct consequence of
Proposition \ref{pr:uni and fact}.

\begin{prop}\label{pr:clockwise}
The following holds
\begin{enumerate}
\item If $c$ does not lie in $\sect(r,r')$, then
\[S_{r'r}=\ccprod_{\ell}S_\ell\]
where $\ell$ ranges over the Stokes rays contained in $\sect(r,r')$, and
$S_\ell$ is placed to the left of $S_{\ell'}$ if $\ell$ is contained in $\sect
(\ell',r')$.
\item If $c$ lies in $\sect(r,r')$, then
 \[S_{r'r}=\ccprod_{\ell}S_\ell\cdot e^{2\pi\iota[B]}\cdot \ccprod_{\ell}S_\ell\cdot e^{-2\pi\iota[B]}\]
where the leftmost product ranges over the Stokes rays contained in
$\sect(c,r')$, and the rightmost one over those contained in $\sect(r,c)$
except for $c$ is the latter is a Stokes ray.
\end{enumerate}
\end{prop}

\subsection{Stokes matrices}\label{ss:Stokes mat class}

Let $r$ be a ray such that both $\pm r$ are admissible, and distinct from
the log cut $\cut$. Assume further that $\cut$ lies in the cone $\sect(-r,r)$.
\comment{This last condition is only necessary so that the monodromy
relation comes out looking nicer.}
By definition, the {\it Stokes matrices} $S^r_\pm$ are the unipotent elements
of $G$ defined by
\[S^r_+=S_{-r\,r}\aand S^r_-=S_{r\,-r}\]

The pair $(A,r)$ determines a partition $\sfPhi=\sfPhi_+\sqcup\sfPhi_-$ of
the root system given by $\sfPhi_\pm=\{\alpha\in\sfPhi|\,\alpha(A)\in\sect
(\pm r,\mp r)\}$. By Proposition \ref{pr:uni and fact}, the Stokes matrices
$S^r_+,S^r_-$ lie in $N_+,N_-$ respectively, where $N_\pm=N_\pm(A,r)
\subset G$ is the unipotent subgroup with Lie algebra $\n_\pm=\bigoplus
_{\alpha\in\sfPhi_\pm}\g_\alpha$. Moreover, if $A$ is fixed, the Stokes
matrices $S^r_\pm$ (and the subgroups $N_\pm$) are locally constant
in $r$, so long as $\pm r$ do not cross a Stokes ray or $\cut$.

\subsection{Connection matrix}\label{ss:connection}

Recall that the connection $\nabla$ is said to be {\it non--resonant} at
$z=\infty$ if none of the eigenvalues of $\ad(B)$ are positive integers.
The following is well--know (see, \eg \cite{Wasow} for $G=GL_n(\IC)$).

\begin{lemma}\label{le:nr dkz}
If $\nabla$ is non--resonant, there is a unique holomorphic function
$\chh_\infty:\IP^1\setminus\{0\}\to G$ such that $\chh_\infty(\infty)=1$ and,
for any determination of $\log z$, the function $\csol_\infty=\chh_\infty
\cdot z^{B}$ is a solution of $\nabla\gamma_\infty=0$.
\end{lemma}

Fix a log cut $\cut$ and, for any admissible ray $r$ distinct from $\cut$,
define the {\it connection matrix} $C_r\in G$ by
\[\csol_\infty=\csol_r\cdot C_r\]
where the identity is understood to hold on $r$. By Proposition \ref
{pr:r and r'}, $C_r$ is locally constant \wrt $r$, so long as $r$ does
not cross a Stokes ray or $\cut$.

\subsection{Monodromy relation}

The connection matrix $C_r$ is related to the Stokes matrices $S_
\pm^r$ by the following {\it monodromy relation}.

\begin{prop}\label{pr:monodromy reln}
The following holds
\[ C_r\cdot e^{2\pi\iota B}\cdot C_r^{-1}=S^r_-\cdot e^{2\pi\iota[B]}\cdot S^r_+\]
\end{prop}
\begin{pf}
By definition of $S^r_\pm$, the monodromy of $\csol_r$ around a
positive loop $p_0$ around $0$ based at a point $z_0\in r$
is the \rhs of the stated identity. On the other hand, the monodromy of
$\csol_\infty$ around $p_0$ is $e^{2\pi\iota B}$. Since $\csol
_r=\csol_\infty\cdot C_r^{-1}$, the former monodromy is
conjugate to the latter by $C_r$.
\end{pf}

\Omit{
Let $\g\nr\subset\g$ be the dense open subset consisting of elements
$B$ such that none of the eigenvalues of $\ad(B)$ are positive integers.
Consider the map $C:\g\nr\to G$ given by
mapping $B$ to the connection matrix of $\nabla=d-\left(\Az/z^2+
B/z\right)$. 
\comment{Remove this paragraph? If keep it need to make sure that
make it clear that $A$ is fixed.}
}

\subsection{The Stokes map}\label{ss:stokes map}

Let $N_\pm\subset G$ be the unipotent subgroups corresponding
to $(A,r)$, and $B_\pm=H\ltimes N_\pm\subset G$ the solvable
subgroups with Lie algebras $\b_\pm=\h\ltimes\n_\pm$. Consider
the fibred product
\[B_+\times_H B_-=\{(b_+,b_-)\in B_+\times B_-|\, \pi_+(b_+)\pi_-(b_-)=1\}\]
where $\pi_\pm:B_\pm\to H$ are the quotient maps. Following \cite
{Bo2}, we define the {\it Stokes map} to be the analytic map $\calS_r:
\g\longrightarrow B_+\times_H B_-$ given by
\[B\longrightarrow
\left((S^r_+)^{-1}\cdot e^{-\iota\pi[B]},
S^r_-\cdot e^{\iota\pi[B]}\right)\]

Note that $B_-\times_H B_+$ maps to $G$ via the map 
$\beta:(b_+,b_-)\to b_+\cdot b_-^{-1}$.\comment{$\beta$ is a principal
bundle over its image with structure group the order two elements in
$H$.} Moreover, by Proposition \ref{pr:monodromy reln}, the composition
$\beta\circ\calS_r$ is the map $\g\to G$ given by 
\[B\longrightarrow 
(S^r_+)^{-1}\cdot e^{-2\pi\iota[B]}\cdot (S^r_-)^{-1}=
C_r\cdot e^{-2\pi\iota B}\cdot C_r^{-1}\]

\subsection{Linearisation of $G^*$}\label{ss:linearisation}

Assume now that $G$ is reductive, and fix a symmetric, non--degenerate,
invariant bilinear form $(\cdot,\cdot)$ on $\g$. The pair of opposite Borel
subalgebras $\b_\pm$ of $\g$ then gives rise to a solution $\sfr\in\b_-\otimes
\b_+$ of the classical \YBE given by
\begin{equation}\label{eq:standard cybe}
\sfr = x_i\otimes x^i + \half{1} t_a\otimes t^a
\end{equation}
where $\{x_i\},\{x^i\}$ are bases of $\n_-,\n_+$ respectively which are
dual \wrt $(\cdot,\cdot)$, and $\{t_a\},\{t^a\}$ are dual bases of $\h$.

The element $\sfr$ gives $\g$ the structure of a \qt Lie bialgebra, with
cobracket $\delta:\g\to\g\wedge\g$ given by $\delta(x)=[x\otimes 1+1
\otimes x,\sfr]$. 
The dual Lie bialgebra $(\g^*,\delta^t,[\cdot,\cdot]^t)$ may be identified, as
a Lie algebra, with
\[\b_+\times_\h\b_-=\{(X_+,X_-)\in\b_+\oplus\b_-|\pi_+(X_+)+\pi_-(X_-)=0\}\]
where $\pi_\pm:\b_\pm\to\h$ is the quotient map.
This endows $G^*=B_+\times_H B_-$ with the structure of a Poisson--Lie
group, which is dual to $G$.

Endow now $\g^*$ with its standard Kirillov--Kostant--Souriau Poisson
structure given by
\[ \{f,g\}(x)=\langle[d_x f,d_xg],x\rangle\]
where $d_xh\in T^*_x\g^*=\g$ is the differential of $h$ at $x$, 
and $[\cdot, \cdot]$ is the Lie bracket on $\g$.

Let $\nnu:\g^*\to\g$ be the identification induced by the bilinear form
$(\cdot,\cdot)$, and set $\nnu^\vee=-1/(2\pi\iota)\nnu$. The following
remarkable result is due to Boalch \cite{Bo1,Bo2}.

\begin{thm}\label{th:Boalch}
The map $\calS\circ\nnu^\vee:\g^*\to G^*$ is a Poisson map, and 
generically a local analytic diffeomorphism.
\end{thm}

In particular, $\calS\circ\nnu^\vee$ gives a linearisation of the
\PL structure on $G^*$. We shall give an alternative proof
of the fact that $\calS\circ\nnu^\vee:\g^*\to G^*$ is a Poisson map
in Section \ref{se:alt boalch}.

\section{Stokes phenomena and quantum groups} \label{se:Ug Stokes}%

This section is an exposition of \cite{TL}. We explain in particular
how the dynamical KZ equations give rise to a twist which kills the
KZ associator. Sections \ref{se:filtered V}--\ref{se:filtered A} contain
background material required to do calculus with values in infinite--dimensional
filtered vector spaces and their endomorphisms. Throughout the
paper, $\sfh,\hbar$ are two formal parameters related by $\hbar
=2\pi\iota\sfh$. 
 
\subsection{Filtered vector spaces}\label{se:filtered V}

Let $\V$ be a vector space over a field $\sfk$ endowed with
a decreasing filtration
\[\V=\V_0\supseteq\V_1\supseteq\V_2\supseteq\cdots\]
and $\imath$ the map $\V\rightarrow\ds{\lim_{\longleftarrow}}\,
\V/\V_n$. Recall that $\V$ is said to be {\it separated} if $\imath$
is injective, and {\it complete} if $\imath$ is surjective.

If $\sfk=\IC$, and the quotients $\V/\V_n$ are finite--dimensional, we
shall say that a map $F:X\to\V$, where $X$ is a topological space
(resp. a smooth or complex manifold) is continuous (resp. smooth
or holomorphic) if its truncations $F_n:X\to\V/\V_n$ are. If $\V$ is
separated and complete, giving such an $F$ amounts to giving
continuous (resp. smooth or holomorphic) maps $F_n:X\to\V/\V_n$
such that $F_n=F_m\mod \V_m/\V_n$, for any $n\geq m$.

\subsection{Filtered endomorphisms}\label{se:filtered E}

Let $\V$ be as in \ref{se:filtered V}, and $\E\subset\End_\sfk(\V)$
the subalgebra defined by
\[\E=\left\{T\in\End_\sfk(\V)|\,T(\V_m)\subseteq\V_m,m\geq 0\right\}\]

Consider the decreasing filtration $\E=\E_0\supseteq\E_1\supseteq
\cdots$ where $\E_n\subset\E$ is the two--sided ideal given by $\E_
n=\{T\in\E|\,\Im(T)\subseteq\V_n\}$. Note that if the quotients $\V/\V
_n$ are finite--dimensional, the same holds for
\[\E/\E_n\cong
\{T\in\End_\sfk(\V/\V_n)|\,T(\V_m/\V_n)\subseteq \V_m/\V_n,\,0\leq m\leq n\}\]
In particular, if $\sfk=\IC$, we may speak of a continuous (resp. smooth,
holomorphic) map with values in $\E$.

\begin{lemma}\hfill
\begin{enumerate}
\item If $\V$ is separated, so is $\E$.
\item If $\V$ is complete, so is $\E$.
\end{enumerate}
\end{lemma}
\begin{pf}
(1) holds because $\bigcap_{n\geq 0}\E_n=\{T\in\E|\,\Im(T)\subseteq
\bigcap_{n\geq 0}\V_n\}$.
(2) Let $T_n\in\E/\E_n$ be such that $T_n=T_m\mod\E/\E_m$ for
any $n\geq m$. It suffices to find $T\in\End_\sfk(\V)$ such that
$T=T_n\mod\E_n$ for any $n\geq 0$, for it then follows that $T\in
\E$. Let $\{v_i\}_{i\in I}$ be a basis of $\V$. For any $i\in I$, $\{T_n
v_i\}$ is a well--defined element of $\lim_n\V/\V_n$. By completeness
of $\V$, there exists $u_i\in\V$ such that $u=T_n v_i\mod\V_n$
for any $n$. Setting $Tv_i=u_i$ gives the required $T$.
\end{pf}

\subsection{Filtered algebras}\label{se:filtered A}

Let $A$ be a $\sfk$--algebra endowed with an increasing algebra
filtration $\sfk=A_0\subseteq A_1\subseteq\cdots$, and $A\fml^o$
the (completed Rees) algebra given by
\[A\fml^o=
\{\sum_{k\geq 0}a_k\hbar^k\in A\fml|\,
a_k\in A_{k}\}\]

Endow $A\fml^o$ with the decreasing filtration
\begin{equation}\label{eq:filtration}
A\fml^o_n=A\fml^o\cap\hbar^nA\fml
\end{equation}
\wrt which it is easily seen to be separated and complete.
Note that each $A\fml^o_n$ is a two--sided, $\IC\fml$--ideal in $A\fml^o$,
and that the quotients 
\[A\fml^o/A\fml^o_n\cong
A_{0}\oplus\hbar A_{1}\oplus \cdots\oplus\hbar^{n-1} A_{n-1}\]
are \fd if $A$ is filtered by \fd subspaces.

\subsection{Example}

We shall be interested in the case when $A=\Ug^{\otimes m}$ is a tensor
power of an enveloping algebra, with filtration given by $A_k=(U\g_{\leq
k})^{\otimes m}$. Then,
\begin{equation}\label{eq:for future}
\Ug^{\otimes m}\fml^o=
\UU'\otimes\UU^{\otimes m-1}
\cap
\UU\otimes\UU'\otimes\UU^{\otimes m-2}
\cap
\cdots
\cap
\UU^{\otimes m-1}\otimes\UU'
\end{equation}
where $\UU=\Ug\fml$ and $\UU'=\Ug\fml^o$.
Note that $\Ug^{\otimes m}\fml^o\cap\Ug^{\otimes m}=\sfk$. However, if $x
\in\Ug_{\leq k}$, $i=1,\ldots,m$, and
\[x^{(i)}=1^{\otimes i-1}\otimes x\otimes 1^{\otimes m-i}\in\Ug^{\otimes m}\]
then $\hbar^{k-1}\negthinspace\ad x^{(i)}$ is a derivation of $\Ug^{\otimes
m}\fml^o$, which preserves the filtration $\Ug^{\otimes m}\fml^o_n$.

\subsection{The dynamical KZ equations}

Let now $\g$ be a complex reductive Lie algebra, $\h\subset\g$ a Cartan
subalgebra, and $(\cdot,\cdot)$ an invariant inner product on $\g$. Let
$\sfPhi=\{\alpha\}\subset\h^*$ be the root system of $\g$ relative to $\h$,
choose $x_\alpha\in\g_\alpha$ for any $\alpha\in\sfPhi$ such that $(x_
\alpha,x_{-\alpha})=1$, and set 
\[\Kalpha=x_\alpha x_{-\alpha}+x_{-\alpha}x_\alpha\]

Endow $\A=\Ug^{\otimes 2}\fml^o$ with the filtration $\A_n=\Ug^{\otimes 2}
\fml^o\cap\hbar^n\Ug^{\otimes 2}$ as in \eqref{eq:filtration}, and filter $\E=
\{T\in\End_{\IC}(\A)|\,T(\A_n)\subseteq\A_n\}$ as in \ref{se:filtered E}. Since
the quotients $\A/\A_n$ and $\E/\E_n$ are finite--dimensional, we may
speak of continuous, smooth or holomorphic functions with values in
$\A$ and $\E$.

The dynamical KZ (DKZ) connection is the $\E$--valued
connection on $\IC$ given by
\begin{equation}\label{eq:DKZ}
\nabla\DKKZ=d-\left(\sfh\frac{\Omega}{z}+\ad\muone\right) dz
\end{equation}
where $\mu\in\h$, $\Omega\in\g\otimes\g$ is the invariant tensor corresponding
to $(\cdot,\cdot)$ and, given an element $a\in\A$, we abusively denote by $a$
the corresponding left multiplication operator $L(a)\in\E$.

\subsection{Fundamental solution at $z=0$} \label{ss:Upsilon zero}

\begin{prop}\label{pr:Fuchs 0}
\hfill
\begin{enumerate}
\item
For any $\mu\in\h$, there is a unique holomorphic function
$H_0:\IC\to\A$ such that $H_0(0,\mu)=1$ and, for any
determination of $\log z$, the $\E$--valued function
\[\Upsilon_0(z,\mu)=e^{z\ad\muone}\cdot H_0(z,\mu)\cdot z^{\sfh\Omega}\]
satisfies $\nabla\DKKZ\Upsilon_0=0$.
%
\item $H_0$ and $\Upsilon_0$ are invariant under the diagonal action of $\h$.
\item $H_0$ and $\Upsilon_0$ are holomorphic functions of $\mu$,
and $\Upsilon_0$ satisfies
\[
\left(d_\h-\half{\sfh}\sum_{\alpha\in\Phi_+}\frac{d\alpha}{\alpha}
\Delta(\Kalpha)-z\ad d\muone\right)\Upsilon_0=
\Upsilon_0\left(d_\h-\half{\sfh}\sum_{\alpha\in\Phi_+}\frac{d\alpha}{\alpha}
\Delta(\Kalpha)\right)
\]
\end{enumerate}
\end{prop}

\subsection{Fundamental solutions at $z=\infty$} \label{ss:Upsilon infty}

Let $\halfplane_\pm=\{z\in\IC|\,\Im(z)\gtrless 0\}$.

\begin{thm}\label{th:Stokes infty}\hfill
\begin{enumerate}
\item
For any $\mu\in\hreg^\IR$, there are unique holomorphic functions
$H_\pm:\halfplane_\pm\to\A$ such that $H_\pm(z,\mu)$ tends to $1$
as $z\to\infty$ in any sector of the form $|\arg(z)|\in(\delta,\pi-\delta)$,
$\delta>0$ and, for any determination of $\log z$, the $\E$--valued
function
\[\Upsilon_\pm(z,\mu)=H_\pm(z,\mu)\cdot z^{\sfh\Omega_0}\cdot e^{z\ad\muone}\]
satisfies $\nabla\DKKZ\Upsilon_\pm=0$.
%
\item $H_\pm$ and $\Upsilon_\pm$ are invariant under the diagonal action of $\h$.
\item $H_\pm$ and $\Upsilon_\pm$ are smooth functions of $\mu$,
and $\Upsilon_\pm$ satisfies
\[
\left(d_\h-\half{\sfh}\sum_{\alpha\in\Phi_+}\frac{d\alpha}{\alpha}
\Delta(\Kalpha)-z\ad d\muone\right)\Upsilon_\pm=
\Upsilon_\pm\left(d_\h-\half{\sfh}\sum_{\alpha\in\Phi_+}\frac{d\alpha}{\alpha}
(\onetwo{\Kalpha})\right)
\]
%
\end{enumerate}
\end{thm}

\subsection{Remark}

The PDEs (3) in Proposition \ref{pr:Fuchs 0} and Theorem \ref{th:Stokes infty}
do not take values in $\E$, since left multiplication by $\sfh\Delta(\Kalpha),
\sfh\Kalpha^{(1)}$ and $\sfh\Kalpha^{(2)}$ does not preserve $\A$. Let,
however, $\A\subsetneq\wt{\A}\subset\Ug^{\otimes 2}\fml$ be the Rees
algebra \wrt the laxer filtration $(\Ug^{\otimes 2})_k=\sum_{a+b=2k}\Ug_
{\leq a}\otimes\Ug_{\leq b}$, 
and $\wt{\E}$ the corresponding algebra of endomorphisms. Then,
$\Upsilon_0,\Upsilon_\pm$, and left multiplication by $\sfh\Delta
(\Kalpha),\sfh\Kalpha^{(1)}$ and $\sfh\Kalpha^{(2)}$ all lie in $\wt{\E}$,
and these PDEs should be understood as holding in $\wt{\E}$.

\subsection{$\IZ_2$--equivariance}\label{ss:Z2}

Let $\U\subset\IC$ be an open subset. For any functions $F:\U\to\A$
and $G:\U\to\E$, define $F^\vee:-\,\U\to\A$ and $G^\vee:-\,\U\to\E$ by
\[F^\vee(z)=e^{z\ad(\onetwo{\mu})}(F(-z)^{21})
\quad\text{and}\quad
G^\vee(z)=e^{z\ad(\onetwo{\mu})}\cdot G(-z)^{21}\] 
where $G(-z)^{21}=(1\,2)\cdot G(-z)\cdot(1\,2)$. If $F,G$ are local solutions
of the dynamical KZ equations with values in $\A$ and $\E$ respectively,
then so are $F^\vee,G^\vee$. 

\begin{lemma}\label{le:Z2}
The following holds
\begin{enumerate}
\item For $z\in\halfplane_\pm$,
\[\Upsilon_0^\vee(z)=\Upsilon_0(z)\cdot e^{\mp\pi\iota\sfh\Omega}\]
\item For $z\in\halfplane_\mp$,
\[ \Upsilon_\pm^\vee(z)=\Upsilon_\mp(z)\cdot e^{\pm\pi\iota\sfh\Omega_0}\]
\end{enumerate}
\end{lemma}
\begin{pf}
(1) The uniqueness of the holomorphic part $H_0$ of $\Upsilon_0$
implies that $(e^{z\ad\muone}\cdot H_0)^\vee=e^{z\ad\muone}\cdot
H_0$. It follows that $\Upsilon_0^\vee(z)=H_0(z)\cdot (-z)^{\sfh\Omega}
=\Upsilon_0(z)\cdot e^{\mp\iota\pi\sfh\Omega_0}$ since $\log(-z)=
\log z\mp\iota\pi$, depending on whether $\Im z\gtrless 0$.

(2) Similarly, for $z\in\halfplane_\mp$,
\[\begin{split}
\Upsilon_\pm^\vee(z)
&=
e^{z\ad(\onetwo{\mu})}\cdot (1\,2)\cdot H_\pm(-z)\cdot e^{-z\ad\muone}\cdot (-z)^{\sfh\Omega_0}\cdot (1\,2)\\
&=
H_\pm^\vee(z)\cdot e^{z\ad\muone}\cdot (-z)^{\sfh\Omega_0}
\end{split}
\]
The uniqueness of $H_\pm$ implies that $H_\pm^\vee=H_\mp$,
from which the result follows.
\end{pf}

\subsection{Another $\IZ_2$--equivariance}

Let $\U\subset\IC$ be an open subset. For any functions $F:\U\to\A$
and $G:\U\to\E$, define $\wt{F}:\U\to\A$ and $\wt{G}:\U\to\E$ by
\[\wt{F}(z)=e^{-z\ad(\onetwo{\mu})}(F(z)^{21})
\quad\text{and}\quad
\wt{G}(z)=e^{-z\ad(\onetwo{\mu})}\cdot (1\,2)\cdot G(z)\cdot(1\,2)\]
If $F,G$ are local solutions of the dynamical KZ equations with parameter
$\mu\in\h$, then $\wt{F},\wt{G}$ are solutions of the DKZ equations
with parameter $-\mu$.

\begin{lemma}\label{le:another Z2}
The following holds
\[\wt{\Upsilon}_0(z;\mu)=\Upsilon_0(z;-\mu)
\aand
\wt{\Upsilon}_\pm(z;\mu)=\Upsilon_\pm(z;-\mu)\]
\end{lemma}
\begin{pf}
By definition,
\[\wt{\Upsilon}_0(z;\mu)=
e^{-z\ad(\onetwo{\mu})}\cdot e^{z\ad\mu^{(2)}}\cdot H_0^{21}(z;\mu)\cdot z^{\hbar\Omega}
=
e^{-z\ad\mu^{(1)}}\cdot H_0^{21}(z;\mu)\cdot z^{\hbar\Omega}
\]
which coincides with $\Upsilon_0(z;-\mu)$ by uniqueness. Similarly,
\begin{multline*}
\wt{\Upsilon}_\pm(z;\mu)
=
e^{-z\ad(\onetwo{\mu})}\cdot H_\pm^{21}(z;\mu)\cdot e^{z\ad\mu^{(2)}}\cdot z^{\hbar\Omega_0}\\
=
H_\pm^{21}(z;\mu)\cdot e^{-z\ad\mu^{(1)}}\cdot z^{\hbar\Omega_0}
=
\Upsilon_\pm(z;-\mu)
\end{multline*}
where the second equality uses the fact that $H_\pm$ is of weight zero, and 
the third follows by uniqueness.
\end{pf}

\subsection{Differential twist}\label{ss:differential twist}

Fix henceforth the standard determination of $\log z$ with a cut
along the negative real axis, and let $\Upsilon_0,\Upsilon_\pm$
be the corresponding fundamental solutions of the dynamical KZ
equations given in \ref{ss:Upsilon zero} and \ref{ss:Upsilon infty}
respectively. We shall consider $\Upsilon_0$ and $\Upsilon_\pm$
as (single--valued) holomorphic functions on $\IC\setminus\IR_
{\leq 0}$.


\begin{defn}
The differential twist is the smooth function $J_\pm:\h\reg^\IR\to\Ug
^{\otimes 2}\fml^o$ defined by
\[J_\pm=\Upsilon_0(z)^{-1}\cdot\Upsilon_\pm(z)\] 
where $z\in\IC\setminus\IR_{\leq 0}$.
\end{defn}

\begin{rem}
$J_\pm$ takes in fact values in $\E$. However, the form of
$\Upsilon_0$ and $\Upsilon_\pm$ shows that 
\[J_\pm=
z^{-\sfh\Omega}\cdot H_0(z)^{-1}\cdot
\exp(-z\ad\muone)\left(H_\pm\right)\cdot z^{\sfh\Omega_0}
\]
so that it is a left multiplication operator. We therefore
abusively identify $J_\pm$ and $J_\pm(1^{\otimes 2})$.
\end{rem}

\begin{prop}\label{pr:Jpm}
The following holds
\[J_-=e^{\pi\iota\sfh\Omega}\cdot J_+^{21}\cdot e^{-\pi\iota\sfh\Omega_0}\]
\end{prop}
\begin{pf}
Let $G^\vee(z)=e^{z\ad(\onetwo{\mu})}\cdot G(-z)^{(1\,2)}$
be the involution defined in \ref{ss:Z2}. By definition, $J_+^{21}=(\Upsilon_0^\vee)^{-1}\cdot\Upsilon_+^\vee$,
where the \rhs is evaluated for $\Im z<0$. By Lemma \ref{le:Z2}, this is
equal to $e^{-\pi\iota\sfh\Omega}\cdot\Upsilon_0^{-1}\cdot\Upsilon_-
\cdot e^{\pi\iota\sfh\Omega_0}$.
\end{pf}

\subsection{}

For any $\mu\in\hreg^\IR$, set $\sfPhi_+(\mu)=\{\alpha\in\sfPhi|\alpha(\mu)>0\}$.

\begin{thm}\label{th:J}\hfill
\begin{enumerate}
\item $J_\pm$ kills the KZ associator $\Phi\KKZ\in U\g^{\otimes 3}\fml^o$,
that is
\[\Phi\KKZ\cdot\Delta\otimes\id(J_\pm)\cdot J_\pm\otimes 1
=\id\otimes\Delta(J_\pm)\cdot 1\otimes J_\pm\]
\item $J_\pm
=1^{\otimes 2}+\half{\hbar}\sfj_\pm$ mod $\hbar^2$, where
\[\sfj_\pm=
\mp\Omega_-
+\frac{1}{\pi\iota}\sum_{\alpha\in\sfPhi_+(\mu)}
(\log\alpha+\gamma)\left(\Omega_\alpha+\Omega_{-\alpha}\right)\]
with $\Omega_\alpha=x_\alpha\otimes x_{-\alpha}$, $\Omega
_\pm=\sum_{\alpha\in\sfPhi_+(\mu)}\Omega_{\pm\alpha}$, and $\gamma
=\lim_n(\sum_{k=1}^n\frac{1}{k}-\log(n))$ the Euler--Mascheroni constant.
In particular,
\begin{equation}\label{eq:J alt}
\sfj_\pm-\sfj_\pm^{21}
=
\Omega_\pm-\Omega_\mp
\end{equation}
\item As a function of $\mu\in\hreg^\IR$, $J_\pm$ satisfies
\[d_\h J_\pm=
\half{\sfh}\sum_{\alpha\in\Phi_+(\mu)}\frac{d\alpha}{\alpha}
\left(\Delta(\Kalpha)J_\pm-J_\pm(\onetwo{\Kalpha})\right)\]
\end{enumerate}
\end{thm}

\begin{rem}
Note that the PDE satisfied by $J_\pm$ is independent of the chamber
which $\mu$ lies in since $d\log\alpha=d\log(-\alpha)$ and $\Kalpha=
\mathcal K_{-\alpha}$. Note also that this PDE takes values in $\A$.
Indeed, although neither the left multiplication operator $L(\sfh\Delta(
\Kalpha))$ nor the right multiplication $R(\sfh\Kalpha^{(1)}+\sfh\Kalpha
^{(2)})$ leaves $\A$ inavariant, the fact that $\Delta(\Kalpha)=\Kalpha^{(1)}+\Kalpha
^{(2)}+2(\Omega_\alpha+\Omega_{-\alpha})$ implies that 
\[L(\sfh\Delta(\Kalpha))-
R(\sfh\Kalpha^{(1)}+\sfh\Kalpha^{(2)})
=
2L(\sfh\Omega_\alpha+\sfh\Omega_{-\alpha})
+
\ad(\sfh\Kalpha^{(1)}+\sfh\Kalpha^{(2)})
\]
which preserves $\A$ since $\sfh\Omega_\alpha
\in\A$, and $\ad(\sfh\Kalpha^{(i)})$ leave $\A$ invariant by \ref{se:filtered A}.
\end{rem}

\Omit{
Give 1--jet of solutions at 0 and infty
Relation btw J_21(mu) and J(-mu)}

\subsection{Quantisation of $(\g,\sfr)$}

Fix a chamber $\C\subset\h\reg^\IR$, and set $\sfPhi_+=\sfPhi_+(\mu)$,
$\mu\in\C$. Let
\[\sfr=
\Omega_++\half{1}\Omega_0=
\sum_{\alpha\in\sfPhi_+}x_\alpha\otimes x_{-\alpha}+\half{1}\Omega_0\]
be the Drinfeld--Sklyanin $r$--matrix corresponding to $\C$, and $(\g,\sfr)$
the corresponding \qt Lie bialgebra.\footnote{Note that the $\sfr$ considered
in \ref{ss:linearisation} corresponds to the {\it anti}fundamental chamber.}

Set $R\KKZ=e^{\hbar\Omega/2}$, and let
\[\left(U\g\fml,\Delta_0,R\KKZ,\Phi\KKZ\right)\]
be the \qtqha structure on $\Ug\fml$ underlying the monodromy of the KZ
equations \cite{Dr4}, where $\Delta_0$ is the standard cocommutative
coproduct on $U\g$. If $\mu\in\C$, the differential twist $J_\pm=J_\pm(\mu)$
allows to twist this structure, and yields a \qt Hopf algebra $(U\g\fml,\Delta
_\pm,R_\pm)$, where\footnote{Note that $\Delta_\pm$ and $R_\pm$
depend on the additional choice of $\mu\in\C$. Specifically, if $\mu_0,
\mu_1\in\C$, $p:[0,1]\to\C$ is a path with $p(0)=\mu_0,p(1)=\mu_1$,
and $a_p\in\Ug\fml_0$ is the holonomy of the Casimir connection
along $p$, then
\[\Delta_\pm(x)(\mu_1)=
a_p^{\otimes 2}\cdot\Delta_\pm(a_p^{-1}xa_p)(\mu_0)\cdot(a_p^{\otimes 2})^{-1}
\aand
R_\pm(\mu_1)=a_p^{\otimes 2}\cdot R_\pm(\mu_0)\cdot(a_p^{\otimes 2})^{-1}\]
In particular, the \qt Hopf algebras corresponding to different values
of $\mu\in\C$ are all isomorphic.}
\[\Delta_\pm(x)=J_\pm^{-1}\cdot\Delta_0(x)\cdot J_\pm
\aand
R_\pm=(J_\pm^{-1})^{21}\cdot R\KKZ\cdot J_\pm\]

\begin{thm}\hfill
\begin{enumerate}
\item $(U\g\fml,\Delta_+,R_+)$ is a quantisation of $(\g,\sfr)$.
\item $(U\g\fml,\Delta_-,R_-)$ is a quantisation of $(\g,\sfr^{21})$.
\item $(U\g\fml,\Delta_\pm,R_\pm)$ is isomorphic, as a \qt Hopf
algebra, to the Drinfeld--Jimbo quantum group corresponding to
$\g$.
\end{enumerate}
\end{thm}
\begin{pf}
(1)--(2) By \eqref{eq:J alt}, the coefficient of $\hbar$ in $R_\pm$ is
$\half{1}(\Omega\pm\Omega_+\mp\Omega_-)$, which is equal to
$\sfr$ for $R_+$ and $\sfr^{21}$ for $R_-$.

(3) This follows, for example, from Drinfeld's uniqueness of the
quantisation of $(\g,\sfr)$ \cite{DrICM} given that the Chevalley
involution of $\g$ clearly lifts to $(U\g\fml,\Delta_\pm,R_\pm)$.
\end{pf}

\begin{rem} It follows from (4) of Theorem \ref{th:J} that
\begin{equation}\label{eq:Rpm}
R_-=R\KKZ^0\cdot R_+^{21}\cdot (R\KKZ^0)^{-1}
\end{equation}
\end{rem}

\section{The $R$--matrix as a quantum Stokes matrix} \label{se:R = S}%

\subsection{Quantum Stokes matrices}\label{ss:quant Stokes mat}

Recall that $\halfplane_\pm=\{z\in\IC|\,\Im(z)\gtrless 0\}$. Define the quantum
Stokes matrices $S_\pm\in\Ug^{\otimes 2}\fml^o$ by
\[\Upsilon_+=\Upsilon_-\cdot S_+\aand
\Upsilon_-\cdot e^{\hbar\Omega_0}=\Upsilon_+\cdot S_-\]
where the first identity is understood to hold in $\halfplane_-$ after $\Upsilon
_+$ has been continued across the ray $\IR_{\geq 0}$, and the second
in $\halfplane_+$ after $\Upsilon_-$ has been continued across $\IR_{\leq 0}$.

\begin{prop}\label{pr:Spm}
The following holds
\begin{enumerate}
\item $S_-=
e^{-\iota\pi\sfh\Omega_0}\cdot
S_+^{21}\cdot
e^{\iota\pi\sfh\Omega_0}$.
\item
$J_+^{-1}\cdot e^{2\pi\iota\sfh\Omega}\cdot J_+=
S_+^{-1}\cdot e^{2\pi\iota\sfh\Omega_0}\cdot S_-^{-1}$
\item As functions of $\mu\in\C$, the quantum Stokes matrices $S_\pm$
satisfy 
\[d_\h S_\pm=
\frac{\sfh}{2}\sum_{\alpha\in\sfPhi_+}\frac{d\alpha}{\alpha}
\left[\onetwo{\Kalpha},S_\pm\right]\]
\end{enumerate}
\end{prop}
\begin{pf}
(1) Let $f$ be a holomorphic function on $\halfplane_\pm$, and denote by
$\P_\pm(f)$ the analytic continuation of $f$ to $\halfplane_\mp$ across the
half--axis $\IR_{\gtrless 0}$. By Lemma \ref{le:Z2}, and the definition
of $S_-$,
\[\P_-(\Upsilon_+^\vee)
=\P_-(\Upsilon_-)\cdot e^{\iota\pi\sfh\Omega_0}
=\Upsilon_+\cdot S_-\cdot e^{-\iota\pi\sfh\Omega_0}\]
On the other hand, if $\imath:\IC\to\IC$ is the inversion $z\to -z$,
\[\begin{split}
\P_-(\Upsilon_+^\vee)
&=
e^{z\ad(\onetwo{\mu})}\cdot(1\,2)\cdot \P_-(\Upsilon_+\circ\imath)\cdot (1\,2)\\
&=
e^{z\ad(\onetwo{\mu})}\cdot(1\,2)\cdot \P_+(\Upsilon_+)\circ\imath\cdot (1\,2)\\
&=
e^{z\ad(\onetwo{\mu})}\cdot(1\,2)\cdot \Upsilon_-\circ\imath\cdot S_+\cdot (1\,2)\\
&=
\Upsilon_-^\vee\cdot S^{21}_+\\
&=
\Upsilon_+\cdot e^{-\iota\pi\sfh\Omega_0}\cdot S^{21}_+
\end{split}\]
where the last identity uses Lemma \ref{le:Z2}.

(2) By construction, the monodromy of the fundamental solution $\Upsilon
_0$ around a positively oriented loop $\gamma_0$ around $0$ is $e^{2\pi
\iota\sfh\Omega}$. Let now $\gamma_\infty$ be a clockwise loop around
$\infty$ based at $x_0\in\halfplane_+$. Since such a loop crosses the negative
real axis before the positive one, the monodromy of $\Upsilon_+$ around
$\gamma_+$ is $S_+^{-1}\cdot e^{2\pi\iota\sfh\Omega_0}\cdot S_-^{-1}$.
The result now follows from the fact that $\gamma_\infty$ is homotopic to
$\gamma_0$, and $\Upsilon_+=\Upsilon_0\cdot J_+$.

(3) follows from the PDE satisfies by $\Upsilon_0$ and $\Upsilon_\pm$.
\end{pf}


\subsection{The $R$--matrix as a quantum Stokes matrix}


\begin{thm}
The following holds
\[R_+
=e^{\pi\iota\sfh\Omega_0}\cdot S_-^{-1}
\aand
R_-
=e^{\pi\iota\sfh\Omega_0}\cdot S_+^{-1}
\]
\end{thm}
\begin{pf}
By definition of $S_+$, $\Upsilon_+=\Upsilon_-\cdot S_+$, when both
$\Upsilon_\pm$ are considered as single--valued functions on $\IC
\setminus\IR_{\leq 0}$. On the other hand, by definition of $J_\pm$,
\[\Upsilon_+=\Upsilon_0\cdot J_+=\Upsilon_-\cdot J_-^{-1}\cdot J_+\]
Using Proposition \ref{pr:Jpm} therefore yields
\[S_+=
e^{\iota\pi\sfh\Omega_0}\cdot(J_+^{-1})^{21}\cdot e^{-\iota\pi\sfh\Omega}\cdot J_+=
e^{\iota\pi\sfh\Omega_0}\cdot(R_+^{-1})^{21}\]
where the last equality uses the fact that $R\KKZ=\exp(\pi\iota\sfh\Omega)
=R\KKZ^{21}$. The first stated identity now follows from (1) of Proposition
\ref{pr:Spm}. The second one follows from the first and \eqref{eq:Rpm}.
\end{pf}

\section{Quantum duality principle and semiclassical limits} \label{se:quantum duality}

\subsection{Quantised universal enveloping algebras}

Let $\sfk$ be a field of characteristic zero, and $\UU$ a quantised universal
enveloping algebra (QUE) over $\sfk$, that is a topologically free Hopf algebra
over $\sfk[[\hbar]]$ such that $\UU/\hbar\UU$ is isomorphic to the enveloping
algebra $U\g$ of a Lie algebra $\g$ over $\sfk$. Then, $\UU$ induces a Lie
bialgebra structure on $\g$, with cobracket $\delta:\g\to\g\otimes\g$ given by
\[\delta(x)=\left.\frac{\Delta(\wt{x})-\Delta^{21}(\wt x)}{\hbar}\right|_{\hbar=0}\]
where $\wt{x}\in\UU$ is an arbitrary lift of $x$.

\subsection{The algebra $\UU'$}

Let $\eta:\IC\fml\to\UU$ and $\eps:\UU\to\IC\fml$ be the unit and counit,
respectively. $\UU$ splits as $\Ker(\eps)\oplus\IC\fml\cdot 1$, with projection
onto the first summand given by $\pi=\id-\eta\circ\eps$. Let $\Delta^{(n)}:
\UU\to\UU^{\otimes n}$ be the iterated coproduct recursively defined by
$\Delta^{(0)}=\eps$, $\Delta^{(1)}=\id$, and $\Delta^{(n)}=\Delta\otimes
\id^{\otimes (n-2)}\circ\Delta^{(n-1)}$ for $n\geq 2$.

Following Drinfeld, define the subspace $\UU'\subset\UU$ by \cite{DrICM,
Ga}
\[\UU'=
\left\{x\in\UU\left|\, \pi^{\otimes n}\circ\Delta^{(n)}(x)\in\hbar^n\UU^{\otimes n},\,n\geq 1\right.\right\}\]
The definition of $\UU'$ extends that of the completed Rees algebra of $\Ug$
to an arbitrary QUE. Specifically, the following holds.

\begin{lemma}\label{le:example}
If $\UU=U\g[[\hbar]]$ with undeformed coproduct, then 
$x=\sum_{n\geq 0}x_n\hbar^n$ lies in $\UU'$ if, and only if the filtration
order of $x_n$ in $U\g$ is less than or equal to $n$.
\end{lemma}
\begin{pf}
It is easy to see that, for any $x_1,\ldots,x_k\in\g$
\[\pi^{\otimes n}\circ\Delta^{(n)}(x_1\cdots x_k)=
\sum_{\substack{I_1\sqcup\cdots\sqcup I_n=\{1,\ldots,k\}\\ |I_i|\neq 0}}
x_{I_1}\otimes\cdots\otimes x_{I_n}
\]
where, for any $I=\{i_1,\ldots,i_m\}$, with $i_1<\cdots<i_m$, we set $x_I=
x_{i_1}\cdots x_{i_m}$. In particular, $\pi^{\otimes n}\circ\Delta^{(n)}(x_1
\cdots x_k)=0$ if $n\geq k+1$. This implies that $\hbar^\ell x_1\cdots x_k
\in\UU'$ if, and only if $k\leq\ell$.
\end{pf}

\subsection{Quantum duality principle}

Assume now that $\g$ is finite--dimensional, let $\left(\g^*,\delta^t,[\cdot,
\cdot]^t\right)$ be the dual bialgebra, and $G^*$ the formal \PLg with
Lie algebra $\g$. By definition, the algebra of functions on $G^*$ is the
topological Poisson Hopf algebra given by $\sfk[[G^*]]
=(U\g^*)^*$. The following result is due to Drinfeld.

\begin{thm}[\cite{DrICM,Ga}]\label{th:duality}
$\UU'$ is a topologically free $\sfk[[\hbar]]$--module, and a sub Hopf 
algebra of $\UU$. Its multiplication is commutative mod $\hbar$, and
$\UU'/\hbar\UU'$ is isomorphic, as a local, complete Poisson Hopf
algebra to $\sfk[[G^*]]$.
\end{thm}

\subsection{The isomorphism $\UU'/\hbar\UU'\cong\sfk[[G^*]]$} \label{ss:iso formula}

If $\UU=\Ug[[\hbar]]$ with undeformed coproduct, then $\delta=0$ and
$\g^*$ has trivial bracket. In this case $G^*$ is the (germ at 0 of the) abelian group $\g^*
$ and, by Lemma \ref{le:example},
$\UU'/\hbar\UU'=\wh{\gr{\Ug}}=\sfk[[\g^*]]$, where $\wh{\cdot}$ is the
graded completion.

More generally, the isomorphism $\UU/\hbar\UU\cong\Ug$ induces
a canonical isomorphism
\[i_\Delta:\UU'/\hbar\UU'\longrightarrow\sfk[[G^*]]\]
as follows \cite[Rem. 3.7]{EH}. 
Identify $U\g^*$ as the quotient of the tensor algebra $T\g^*$ endowed
with the standard concatenation product and (cocommutative) shuffle
coproduct, and $(U\g^*)^*$ with a sub Hopf algebra of its dual $(T\g^*)
^*=\wh{T\g}=\prod_{n\geq 0} \g^{\otimes n}$, where
the latter is endowed with the (commutative) shuffle product
and deconcatenation coproduct. Then, the isomorphism $i_\Delta:\UU'/\hbar
\UU'\to\sfk[[G^*]]=(U\g^*)^*\subset\wh{T\g}$ is given by noticing that
if $x\in\UU'$, $\left(\left.\frac{1}{\hbar^n}  \pi^{\otimes n}\circ\Delta^{(n)}
(x)\right)\right |_{\hbar=0}$ lies in $\g^{\otimes n}\subset (\Ug)^{\otimes n}$
for any $n$, and setting
\begin{equation}\label{eq:iso formula}
i_\Delta(x)=
\left\{\left.\frac{\pi^{\otimes n}\circ\Delta^{(n)}(x)}{\hbar^n}\right |_{\hbar=0}\right\}_{n\geq 0}
\in
\prod_{n\geq 0} \g^{\otimes n}
\end{equation}

\subsection{Semiclassical limit}

\comment{What are $G/\g$? At the beginning just algebraic groups,
but pretty soon reductive too.}
If $\UU$ is a QUE which deforms $U\g$, and $A\in\UU\otimes\UU'$,
we denote by $\scl{A}$ the {\it semi--classical limit} of $A$, that is its
class in $\UU\otimes\UU'/\hbar\UU\otimes\UU'$. By Theorem \ref
{th:duality}, $\scl{A}$ lies in $U\g\wh{\otimes}\sfk[[G^*]]$, and is
therefore a (formal) function on $G^*$ with values in $U\g$.
\comment{Mention completed tensor products somewhere}

\section{Semiclassical limit of the dynamical KZ equation}\label{se:scl DKZ}


The goal of this section is to prove that the Stokes data of the ODE
\eqref{eq:dKZ} are the semiclassical limits of the Stokes data of the
dynamical KZ equations \eqref{eq:DKZ}. Technicalities aside, this
stems from the observation that if $\Upsilon$ is a solution of
\[\frac{d\Upsilon}{dz}=
\left(\ad\muone+\sfh\frac{\Omega}{z}\right)\Upsilon\]
with values in $\UU\otimes\UU'$, the semiclassical limit $\gamma$
of $\Upsilon$, as a formal function of $\lambda\in\g^*$ with values
in $\Ug$, satisfies
\[\frac{d\gamma}{dz}=
\left(\ad\mu+\frac{\nnu(\lambda)}{2\pi\iota z}\right)\gamma\]
where $\nnu(\lambda)=\id\otimes\lambda(\Omega)$ which, after
the change of variable $z\to 1/z$, and the replacement $\ad\mu
\to -\Az, \nnu(\lambda)\to -2\pi\iota B$ is precisely the equation \eqref
{eq:dKZ}.\footnote{The appearance of the factor $2\pi\iota$
is due to the fact that the identification $\UU'/\hbar\UU'\cong\wh{S\g}$
is given by mapping $x\in\g$ to $\hbar x=2\pi\iota\sfh x\in\UU'$.}

\subsection{Formal Taylor series groups}  

Let $G$ be an affine algebraic group over $\IC$. The ring of
regular functions $\IC[G]$ is a Hopf algebra, with coproduct
$\Delta f(g_1,g_2)=f(g_1g_2)$, counit $\eps(f)=f(1)$, and
antipode $Sf(g)=f(g^{-1})$.

If $(\RA,m_{\RA},1_{\RA})$ is a commutative, unital $\IC
$--algebra, the $\RA$--points of $G$ are, by definition, the
set of $\IC$--algebra morphisms $G(\RA)=\Alg_\IC(\IC[G],\RA)$.
$G(\RA)$ is a group, with multiplication $\phi\cdot\psi=m_
{\RA}\circ\phi\otimes\psi\circ\Delta$, unit $1_{\RA}\circ
\epsilon$, and inverse $\phi^{-1}=\phi\circ S$. Let $\mm\subset
\RA$ be a maximal ideal, and denote by $G(\RA)_\mm\subset
G(\RA)$ the normal subgroup consisting of maps $\gamma:\Spec
R\to G(\IC)$ such that $\gamma(\mm)=1$, that is
\[G(\RA)_\mm=
\left\{\varphi\in\Alg_\IC(\IC[G],\RA)|\,\varphi(I)\subset\mm\right\}\] 
where $I=\Ker\eps$ is the augmentation ideal. We shall need the
following elementary

\begin{lemma}
If $\RA$ is a complete local ring with unique maximal ideal $\mm$,
then $G(\RA)_\mm$ may be identified with the set of grouplike
elements of the topological Hopf algebra
$$\Ug\wh{\otimes}\RA=\lim_p \Ug\otimes\RA/\mm^p$$
\end{lemma}
\begin{pf}
Let $\IC[[G]]=\lim\IC[G]/I^n$ be the completion of $\IC[G]$ at the
identity, and identify $\Ug$, as a Hopf algebra, with the continuous
dual
\[\IC[[G]]^*=\left\{\varphi\in\Hom_\IC(\IC[G],\IC)|\,\varphi(I^n)=0, n\gg 0\right\}\]
If $\mm^p=0$ for some $p$, and $\phi\in G(\RA)_\mm$, $\phi$ vanishes
on $I^p$ and therefore lies in $\left(\IC[G]/I^p\right)^*\otimes\RA\subset
\IC[[G]]^*\otimes\RA$. In general, $\mm$ is of finite order in $\RA/\mm^p$
for any $p\geq 1$, so that $G(\RA)_\mm=\lim_p G(\RA/\mm^p)_\mm$
embeds into $\lim_p\Ug\otimes\RA/\mm^p$.
\end{pf}

We shall be interested below in the case when $\RA=\IC[[V]]$ is
the completion of the algebra of regular functions on the vector
space $V=\g$ or $V=\g^*$ at $0$. We denote in this case $G(\RA),
G(\RA)_\mm$ and $\Ug\wh{\otimes}\RA$ by $G[[V]],G[[V]]_+$
and $\Ug[[V]]$ respectively. As algebraic groups over $\IC$,
$G[[V]]$ and $G[[V]]_+$ are the inverse limits
\[ G[[V]]=\lim_{\longleftarrow} G[[V]]^{(m)} 
\aand
G[[V]]_+=\lim_{\longleftarrow} G[[V]]^{(m)}_+\] 
where $G[[V]]^{(m)}=G(\IC[[V]]/I^m)$, respectively, and $G[[V]]_+$
is prounipotent.

\subsection{Semiclassical limit of canonical solutions of the DKZ equations}

Consider the ODE
\begin{align}
\frac{d\gamma}{dz}
&=
\left(\frac{A}{z^2}+\frac{B}{z}\right)\gamma
\label{eq:class}\\
\intertext{and the dynamical KZ equation}
\frac{d\Upsilon}{dz}
&=
\left(\ad\muone+\sfh\frac{\Omega}{z}\right)\Upsilon
\label{eq:quant}
\end{align}
where $A,\mu\in\hreg$, and $B\in\g$.

Fix throughout the standard determination of the logarithm, with a cut along
$\IR_{<0}$. The following result shows that the semiclassical limits of the
canonical fundamental solutions of \eqref{eq:quant} at $z=0,\infty$ are the
canonical fundamental solutions of \eqref{eq:class} at $z=\infty,0$, after
the change of variable $z\to 1/z$. 

\begin{prop}\label{pr:scl DKZ}
Let $\nu:\g^*\to\g$ be the isomorphism given by $\lambda\to\lambda\otimes
\id(\Omega)$, and set $\nu^\vee=-\nu/2\pi\iota$.
\begin{enumerate}
\item 
Let $\csol_\infty$ be the canonical solution of \eqref{eq:class} near
$z=\infty$, and write $$\csol_\infty=e^{-\Az/z}\cdot h_\infty\cdot z^B$$
where $h_\infty:\IP^1\setminus 0\to G$ is such that $h_\infty(\infty)=1$.
Regard $h_\infty$ as a holomorphic function of $B\in\g\nr$ such that
$\left.h_\infty(z)\right|_{B=0}\equiv 1$, and let
\[\whh_\infty:\IP^1\setminus 0\longrightarrow G[[\g]]_+\]
be its formal Taylor series at $B=0$.

Let $\qsol_0=e^{z\ad\muone}\cdot H_0\cdot z^{\sfh\Omega}$ be
the canonical solution of \eqref{eq:quant} near $z=0$. Then, the
semiclassical limit of $H_0$ takes values in $G[[\g^*]]
_+\subset\Ug[[\g^*]]$. Moreover, if $\mu=-\Az$, then
\[\scl{H_0(z)}(\lambda)=\whh_\infty(1/z;\nu^\vee(\lambda))\]

\item 
Assume now that $\Az\in\hreg^\IR$. Let $$\csol_\pm=h_\pm\cdot
e^{-\Az/z}\cdot z^{[B]}:\halfplane_\pm\to G$$ be the canonical solution
of \eqref{eq:class} at $z=0$ corresponding to the
half--plane $\halfplane_\pm=\{z\in\IC|\,\Im(z)\gtrless 0\}$. Regard $h_\pm$ as a holomorphic
function of $B\in\g$ such that $\left.h_\pm(z)\right|_{B=0}\equiv 1$,
and let
\[\whh_\pm:\halfplane_\pm\longrightarrow G[[\g]]_+\]
be its formal Taylor series at $B=0$.

Let $\qsol_\pm=H_\pm\cdot e^{z\ad\muone}\cdot z^{\sfh\Omega_0}$,
be the canonical solution of \eqref{eq:quant} at $z=\infty$ corresponding
to the half--planes $\halfplane_\pm$. Then, the semiclassical limit of
$H_\pm$ takes values in $G[[\g^*]]_+\subset\Ug[[\g^*]]$. Moreover,
if $\mu=-\Az$, then
\[\scl{H_\pm(z)}(\lambda)=\whh_\mp(1/z;\nu^\vee(\lambda))\]
\Omit{
Let $H_\pm:\halfplane_\pm\to\A$ be the holomorphic component of the
canonical solution of the DKZ equations at $z=\infty$ corresponding
to the half--planes $\halfplane_\pm$. Then, the semiclassical limit $\scl{H_
\pm}$ of $H_\pm$ takes values in $G[[\g^*]]_+\subset\Ug[[\g^*]]$.
Moreover, if $\mu=-\Az$, then
\[\scl{H_\pm(z)}(B,\mu)=\whh_\mp(1/z;-B,\Az)\]
}
\end{enumerate}
\end{prop}

\begin{pf}
(1)
By definition, $H_0$ is a solution of 
\[\frac{dH_0}{dz}=\frac{\sfh}{z}\left(\ell(e^{-z\ad\muone}(\Omega))-\rho(\Omega)\right)H_0\]
where $\ell,\rho$ denote left and right multiplication respectively. Thus,
as en element of $\UU\otimes\UU'/\hbar\UU\otimes\UU'=\Ug[[\g^*]]$,
the semiclassical limit $h_0$ of $H_0$ satisfies
\[\frac{dh_0}{dz}=\frac{1}{2\pi\iota z}\left(\ell(e^{-z\ad\mu}(\nnu)-\rho(\nnu)\right)h_0\]
together with the initial condition
$h_0(0)=1$. We claim that $h_0$ takes values in $G[[\g^*]]_+\subset
\Ug[[\g^*]]$. Indeed, both $\Delta\otimes\id(h_0)$ and $h_0^{13}h_0^{23}$
satisfy
\[\frac{dh}{dz}=\frac{1}{2\pi\iota z}\left(\ell(e^{-z\ad\mu}(\nnu^1+\nnu^2))-\rho(\nnu^1+\nnu^2)\right)h\]
and the result follows by uniqueness. The claimed equality now follows
from the uniqueness statement of Lemma \ref{le:nr dkz}, applied to the affine
algebraic groups $G[\IC[\g^*]/I^m]$, $m\geq 1$.

(2) is proved similarly.
\end{pf}

\subsection{Semiclassical limit of the differential twist}\label{ss:scl J} 

\begin{thm}\label{th:scl J}
Assume that $\Az\in\hreg^\IR$, and let $C_\pm=\csol_\pm^{-1}\cdot\csol
_\infty$ be the connection matrix of \eqref{eq:class} (see \ref{ss:connection}).
Regard $C_\pm$ as a $G$--valued holomorphic function of $B\in\g\nr$
such that $\left.C_\pm\right|_{B=0}=1$, and let $\wh{C}_\pm\in G[[\g]]_
+$ be its formal Taylor series at $B=0$.

Let $J_\pm=\Upsilon_0^{-1}\cdot\Upsilon_\pm$ be the differential
twist defined in \ref{ss:differential twist}. Then, if $\mu=-A$, the
semiclassical limit of $J_\pm$ is given by
\[\scl{J_\pm}(\lambda)=\wh{C}_\mp(\nu^\vee(\lambda))^{-1}\]
\end{thm}
\begin{pf}
By definition,
$J_\pm=
z^{-\sfh\Omega}\cdot H_0(z)^{-1}\cdot
\exp(-z\ad\muone)\left(H_\pm\right)\cdot z^{\sfh\Omega_0}
$, 
where $z\in\halfplane_\pm$. By Proposition \ref{pr:scl DKZ},
\[\scl{J_\pm}=
z^{-\nnu/2\pi\iota}\cdot \wh{h}_\infty(1/z;-\nnu/2\pi\iota)^{-1}\cdot \exp(z\ad(\Az))
\left(\wh{h}_\mp(1/z;-\nnu/2\pi\iota)\right)\cdot z^{[\nnu]/2\pi\iota}
\]
On the other hand,
\[C_\pm(B)=
w^{-[B]}\cdot e^{\Az/w}\cdot h_\pm(w)^{-1}\cdot
e^{-\Az/w}\cdot h_\infty(w)\cdot w^B\]
where $w\in\halfplane_\pm$.
\end{pf}

\subsection{Semiclassical limit of the quantum Stokes matrices}
\label{ss:sc Stokes mat}

\begin{thm}
Let $A\in\hreg^\IR$, and $S_\pm$ the Stokes matrices of
the ODE \eqref{eq:class} relative to the ray $-\iota\IR_{>0}$
(see \ref{ss:Stokes mat class}). Regard $S_\pm$ as a $G$--valued
holomorphic function of $B\in\g$ such that $\left.S_\pm\right|
_{B=0}=1$, and let $\wh{S}_\pm\in G[[\g]]_+$ be its formal
Taylor series at $0$. 

Let $\mu\in\hreg^\IR$, and $S^\hbar_\pm$ the Stokes matrices
of the dynamical KZ equation \eqref{eq:quant} (see \ref {ss:quant Stokes mat}).
Then, $S^\hbar_\pm$ take values in $\UU\wh{\otimes}\UU'$,
and its semi--classical limit in $G[[\g^*]]_+\subset\Ug[[\g^*]]$.
Moreover, if $\mu=-A$, then
\[\scl{S^\hbar_\pm}(\lambda)=\wh{S}_\pm(\nu^\vee(\lambda))\]

%
\end{thm}
\begin{pf}
Let $\qsol_\pm=H_\pm\cdot e^{z\ad\muone}\cdot z^{\sfh\Omega_0}$
be the canonical solutions of the DKZ equations corresponding to
the halfplanes $\halfplane_\pm$, and $\wt{\qsol}_+=\wt{H}_+\cdot e^{z\ad
\muone}\cdot z^{\sfh\Omega_0}$ the analytic continuation of $\qsol_+$
across $\IR_{>0}$. By definition,
\[S_+^\hbar=
\qsol_-^{-1}\cdot\wt{\qsol}_+=
z^{-\sfh\Omega_0}\cdot \exp(-z\ad\muone)\left(H_-^{-1}\cdot
\wt{H}_+\right)\cdot z^{\sfh\Omega_0}\]
for $z\in\halfplane_-$. By Proposition \ref{pr:scl DKZ},
\[\scl{S^\hbar_+}=z^{-[\nnu]/2\pi\iota}\cdot\exp(-z\ad(\mu))\left(\wh{h}_+(1/z;-\nnu/2\pi\iota)^{-1}\cdot\wt{\wh{h}}_-(1/z;-\nnu/2\pi\iota)\right)z^{[\nnu]/2\pi\iota}\]
where $\wt{\wh{h}}_-$ is the analytic continuation of $\wh{h}_-$
across $\IR_{>0}$.

On the other hand, if $\gamma_\pm(w)=h_\pm(w)\cdot e^{-A/w}\cdot w^{[B]}$
are the canonical solutions of \eqref{eq:class} corresponding to $w\in\IH_\pm$,
and $\wt{\gamma}_-$ is the analytic continuation of $\gamma_-$ across
$\IR_+$ then, by definition
\[S_+=
\csol_+^{-1}\cdot\wt{\csol}_-=
w^{-[B]}\cdot e^{A/w}\cdot h_+(w)^{-1}\cdot \wt{h}_-(w)\cdot e^{-A/w}\cdot w^{[B]}
\]
The Taylor series of $S_+$ at $B=0$ therefore coincides with $\scl{S^\hbar_+}$
provided $A=-\mu$, $w=1/z$, and $B$ is replaced by $-\nnu(\lambda)/2\pi\iota$. The proof that $\scl
{S^\hbar_-}=\wh{S}_-(-\nnu/2\pi\iota)$ is identical.
\end{pf}

\section{Formal linearisation via quantisation} \label{se:EEM}%

\subsection{}
Let $(\p,\sfr)$ be a \fd \qt Lie bialgebra over a field $\sfk$ of characteristic
zero. Thus, $\p$ is a Lie algebra, $\sfr\in\p\otimes\p$ satisfies the classical
\YBE (CYBE)
\[[\sfr_{12},\sfr_{23}+\sfr_{13}]+[\sfr_{13},\sfr_{23}]=0\]
and is such that $\Omega=\sfr+\sfr^{21}$ is invariant under $\p$. In particular,
$\p$ is a Lie bialgebra with cobracket $\delta:\p\to\p\wedge\p$ given by
$\delta(x)=[x\otimes 1+1\otimes x,\sfr]$. 

Let $\p^*$ be the dual Lie bialgebra to $\p$, and $P,P^*$ the formal
\PL groups with Lie algebras $\p,\p^*$. The CYBE imply that the maps
$\ell,\rho:\p^*\to\p$ given by
\[\ell(\lambda)=\lambda\otimes\id(\sfr)\aand \rho(\lambda)
=-\id\otimes\lambda(\sfr)\]
are morphisms of Lie algebras. We denote the corresponding morphisms
of formal groups $P^*\to P$ by $\LL$ and $\RR$ respectively, and by $\beta:P^*
\to P$ the {\it big cell map}
\[g^* \longrightarrow \LL(g^*)\cdot\RR(g^*)^{-1}\]
The differential of $\beta$ at $1$ is $\ell-\rho:\lambda\to\lambda\otimes\id(\Omega)=
:\nnu(\lambda)$. In particular, $\beta$ is an isomorphism of formal manifolds
if $\sfr$ is {\it non--degenerate}, that is such that $\nnu:\p^*\to\p$ is an isomorphism.

\subsection{}

Set $\UU=U\p[[\hbar]]$ and let $\Phi\in\UU^{\otimes 3}$ be an associator,
that is an element satisfying $\Phi\in 1+\frac{\hbar^2}{24}[\Omega_{12},\Omega
_{23}]+\hbar^3\UU^{\otimes 3}$, and such that $(\UU,\Delta_0,e^{\hbar\Omega/2},
\Phi)$ is a \qt quasi--Hopf algebra. Let $J\in 1+\frac{\hbar}{2}\sfj+\hbar^2\UU^
{\otimes 2}$ be a twist such that $\sfj-\sfj^{21}=\sfr-\sfr^{21}$, and the following
twist equation holds
\begin{equation}\label{eq:twist eq}
\Phi\cdot J_{12,3}\cdot J_{1,2}=J_{1,23}\cdot J_{2,3}
\end{equation}
Then, 
$\UU_J=\left(\UU,\Delta_J=J^{-1}\Delta_0(\cdot)J,R_J=J^{-1}_{21}e^{\hbar\Omega/2}J\right)$
is a \qt Hopf algebra, which is a quantisation of $(\p,\sfr)$. By Theorem \ref{th:duality},
$(\UU_J)'$ is therefore a quantisation of the Poisson algebra $\sfk[[P^*]]$.

\subsection{}

Assume that the twist $J$ is {\it admissible}, that is such that $\hbar\log(J)\in(\UU')
^{\otimes 2}$. The following linearisation result is due to Enriquez--Halbout \cite
[Prop. 4.2]{EH}.
\begin{prop}
The subalgebras $\UU'$ and $(\UU_J)'$ of $\UU$ coincide. Their equality therefore
induces a formal Poisson isomorphism $\poiss:\p^*\to P^*$ given by the composition
\begin{equation}\label{eq:EH}
\sfk[[P^*]]\cong(\UU_J)'/\hbar(\UU_J)' = \UU'/\hbar\UU'\cong\sfk[[\p^*]]
\end{equation}
where the first and last isomorphisms are given by \eqref{eq:iso formula}.
\end{prop}

The explicit form of the isomorphism $\poiss$ is given by the following
result of \EEM \cite[\S 3.3.2]{EEM}.
\begin{thm}\label{th:EEM}
Assume further that $\Phi=\Psi(\hbar\Omega_{12},\hbar\Omega_{23})$
where $\Psi$ is a Lie associator, and that $J$ lies in $\UU\otimes\UU'
\cap\UU'\otimes\UU$. Then,
\begin{enumerate}
\item The \sscl $\jmath=J\mod\hbar\,\UU\otimes\UU'$ lies in $P[[\p^*]]_+
\subset\Ug\wh{\otimes}\sfk[[\p^*]]$, that is, is a formal map $\p^*\to P$.
\item The composition of the Poisson isomorphism $\poiss:\p^*\to P^*$ 
with the big cell map $\beta:P^*\to P$ is the map $e_\jmath:\p^*\to P$ defined by
\[e_\jmath(\lambda)=\jmath(\lambda)^{-1}\cdot e^{\nnu(\lambda)}\cdot\jmath(\lambda)\]
\end{enumerate}
In particular, if $\sfr$ is non--degenerate, the map $\beta^{-1}\circ e
_\jmath:\p^*\to P^*$ is an isomorphism of formal Poisson manifolds.
\end{thm}
%
\begin{pf}
We outline the proof for the reader's convenience. 
By assumption, $R=J_{21}^{-1}\cdot e^{\hbar\Omega/2}\cdot J$ lies in
$\UU\otimes\UU'=\UU_J\otimes(\UU_J)'$, and similarly $R_{21}\in\UU
_J\otimes(\UU_J)'$.

Consider now the identity
\begin{equation}\label{eq:basic id}
R_{21}\cdot R=J^{-1}\cdot e^{\hbar\Omega}\cdot J
\end{equation}
Let $b\in\UU_J\otimes(\UU_J)'/\hbar\UU_J\otimes(\UU_J)'\cong U\p[[P^*]]$ be
the \sscl of the left--hand side, and $a\in\UU\otimes\UU'/\hbar\UU\otimes
\UU'\cong U\p[[\p^*]]$ that of the right--hand side. Clearly, $b\circ\poiss=a$. It
therefore suffices to show that $b=\beta$ and $a=e_\jmath$.

The identity $\Delta_J\otimes\id(R)=R_{13}\cdot R_{23}$ implies that the
\sscl $R'$ of $R$ lies in $P[[P^*]]_+$, and $\id\otimes\Delta_J(R)=R_{13}
\cdot R_{12}$ that $R'$ is an antihomomorphism $P^*\to P$. Its differential
at the identity is readily seen to be the map $\p^*\to\p$ given by $\lambda
\to\id\otimes\lambda(\sfr)$, so that $R'(g^*)=\RR(g^*)^{-1}$. Similarly, the 
\sscl of $R^{21}$ is the homomorphism $\LL:P^*\to P$, and it follows that
$b=\beta$.\comment{Notice an interesting peculiarity. The big cell map
$\beta:P^*\to P$ is not a Poisson map (afaik) and therefore one does
not expect to quantise it to an algebra homomorphism 
$$\beta^*_\hbar:\sfk[[P]]_\hbar=(\UU_J)^*\longrightarrow(\UU_J)'=\sfk[[P^*]]_\hbar$$ or equivalently an element $B^*_\hbar\in
\UU_J\otimes(\UU_J)'$ satisfying $\Delta_J\otimes\id(B^*_\hbar)=(B^*_\hbar)
_{13}(B^*_\hbar)_{23}$. However, $\beta$ {\bf is} being quantised as the element
$RR_{21}$ of $\UU_J\otimes(\UU_J)'$. Clarify this.}

Since the \sscl of $e^{\hbar\Omega}$ is $e^\nu\in P[[\p^*]]_+$, we have
$a=e_\jmath$ and there remains to prove that $\jmath$ lies in $P[[\p^*]]
_+$, that is satisfies $\Delta_0\otimes\id(\jmath)=\jmath_{1,3}\cdot\jmath
_{2,3}$. This is a consequence of the reduction of the twist equation \eqref
{eq:twist eq} mod $\hbar\UU\otimes\UU\otimes\UU'$, as follows. 
Note first that since $J\in1+\hbar\UU\otimes\UU$, $J_{1,2}\in 1+\hbar\UU
\otimes\UU\otimes\UU'$. Next, it is easy to see that for any $x\in\UU'$,
$\Delta_0(x)\in 1\otimes x+\hbar\UU\otimes\UU'$, hence $J_{1,23}\in
J_{1,3}+\hbar\UU\otimes\UU\otimes\UU'$. Finally, $\hbar\Omega_{12}\in
\hbar\UU\otimes\UU\otimes\UU'$, hence $\Phi=\Psi(\hbar\Omega_{12},
\hbar\Omega_{23})=\Psi(0,\hbar\Omega_{23})=1$ mod $\hbar\UU\otimes
\UU\otimes\UU'$.
\end{pf}

\section{Analytic linearisation via Stokes data} \label{se:alt boalch}%

Let $G$ be a complex reductive group, and $B_\pm\subset
G$ a pair of opposite Borel subgroups intersecting along the
maximal torus $H$. Let $\g,\b_\pm,\h$ be the Lie algebras of
$G,B_\pm$, and $H$ respectively, $\sfPhi\subset\h^*$ the corresponding
root system, and $\sfPhi_\pm\subset\sfPhi$ the set of positive
and negative roots, so that $\b_\pm=\h\bigoplus_{\alpha\in\sfPhi
_\pm}\g_\alpha$. Fix an invariant inner product $(\cdot,\cdot)$
on $\g$, and let $\sfr\in\b_+\otimes\b_-$ be the corresponding
canonical element (see \eqref{eq:standard cybe}). Then, $(\g,
\sfr)$ is a quasitriangular Lie bialgebra, and $G$ and $G^*=B
_-\times_H B_+$ are dual Poisson--Lie groups. Moreover, the
homomorphisms $L,R:G^*\to G$ defined in \ref{ss:EEM}
correspond to the first and second projection, respectively.

Let $A\in\h\reg$, and consider the connection
\[\nabla=d-\left(\frac{A}{z^2}+\frac{B}{z}\right)dz\]
Set $\h^\IR=\{t\in\h|\alpha(t)\in\IR,\alpha\in\sfPhi\}$, and
let $\C=\{t\in\h|\,\alpha(t)>0,\,\alpha\in\sfPhi_+\}\subset\h^\IR\reg$
be the fundamental chamber corresponding to $\sfPhi_+$.  
Note that the rays $\pm\iota\IR_{>0}$ are admissible if $A\in\h
\reg^\IR+\iota\h^\IR\subset\h\reg$. Moreover, by \ref{ss:Stokes mat class},
the Stokes matrices $S_\pm$ corresponding to $r=-\iota\IR_{>0}$
lie in $B_\pm(A,r)=B_\mp$ if $A\in -\C+\iota\h^\IR$. Let
\[\calS:\g\to G^*
\qquad\qquad
B\to\left(S_+^{-1} \cdot e^{-\iota\pi[B]}, S_-\cdot e^{\iota\pi[B]}\right)\]
be the Stokes map defined in \ref{ss:stokes map}.

Let $\nu:\g^*\to\g$ be the identification determined by $(\cdot,\cdot)$,
and set $\nnuc=-\nu/2\pi\iota$.

\begin{thm}
If $A\in -\C$, the map $\calS\circ\nnuc:\g^*\to G^*$ is a Poisson map. 
\end{thm}
\begin{pf}
Since $\calS\circ\nnuc$ is complex analytic, it is sufficient to prove that
its formal Taylor series at $0$ is a Poisson map.
%

Set $\mu=-A\in\C$, and let $J_+=J_+(\mu)$ the differential twist defined
in \ref{ss:differential twist}. By Theorem \ref{th:J}, $J_+\in 1+\half{\hbar}
\sfj_++\hbar^2\UU^{\otimes 2}$, where $\sfj_+-\sfj^{21}_+=\sfr-\sfr^{21}$,
and $J_+$ kills the KZ associator $\Phi\KKZ$.

Write $\Omega=\Omega_0+\sum_{\alpha\in\sfPhi}\Omega_\alpha$, 
where $\Omega_0=\sum_i t_i\otimes t^i$, with $\{t_i\},\{t^i\}$ dual
bases of $\h$ \wrt $(\cdot,\cdot)$, and $\Omega_\alpha=x_\alpha
\otimes x_{-\alpha}$, with $x_{\pm\alpha}\in\g_{\pm\alpha}$ such
that $(x_\alpha,x_{-\alpha})=1$. Then, one can show that $\log J_+$
is a Lie series in the variables $\hbar\Omega_0,\hbar\Omega_\alpha$.
Since the subspace $\adm_n=\{x\in\UU^{\otimes n}|\,\hbar x\in(\UU')
^{\otimes n}\}$ is a Lie algebra for any $n\geq 1$, and $\hbar\Omega
_0,\hbar\Omega_\alpha\in\adm_2$, it follows that $\log J_+\in\adm_2$.

Since $J_+$ lies in $\UU'\otimes\UU\cap\UU\otimes\UU'$ by \ref{ss:differential twist},
we may apply Theorem \ref{th:EEM} to the pair $(\Phi\KKZ,J_+)$. Let
$\jmath_+\in\UU\otimes\UU'/\UU\otimes\UU'=G[[\g^*]]_+$ be the \sscl
of $J_+$, and $e_{\jmath_+}\in G[[\g^*]]_+$ the map $\lambda\to\jmath
_+(\lambda)^{-1}\cdot e^{\nu(\lambda)}\cdot\jmath_+(\lambda)$. By
Theorem \ref{th:scl J}
\[\begin{split}
e_{\jmath_+}(\lambda)
&=
\wh{C}_-(-\nu(\lambda)/2\pi\iota;-\mu)\cdot e^{\nu(\lambda)}\cdot \wh{C}_-(-\nu(\lambda)/2\pi\iota;-\mu)^{-1}\\
&=
\left(
\wh{C}_-(\nu^\vee(\lambda);A)\cdot e^{2\pi\iota\nu^\vee(\lambda)}\cdot \wh{C}_-(\nnu^\vee(\lambda);A)^{-1}
\right)^{-1}\\
&=
\left(
\wh{S}_-(\nu^\vee(\lambda);A)\cdot e^{2\pi\iota[\nu^\vee(\lambda)]}\cdot\wh{S}_+(\nnu^\vee(\lambda);A)
\right)^{-1}\\
&=
\left(\wh{S}_+(\nnu^\vee(\lambda);A)^{-1}\cdot e^{-\pi\iota[\nu^\vee(\lambda)]}\right)
\cdot
\left(\wh{S}_-(\nnu^\vee(\lambda);A)\cdot e^{\pi\iota[\nu^\vee(\lambda)]}\right)^{-1}\\
&=
L(\wh{\calS}(\nu^\vee(\lambda);A))\cdot R(\wh{\calS}(\nu^\vee(\lambda);A))^{-1}
\end{split}\]
where the 
the third equality follows from
the monodromy relation (Proposition \ref{pr:monodromy reln}), and the
fifth from the definition of the Stokes map, as well as the assumption that
$A\in-\C$, so that $S_\pm(B;A)\in N_\pm(A,r)=N_\mp$.

It follows that the composition $\beta^{-1}\circ e_{\jmath_+}$ is equal to
$\wh{\calS}\circ\nu^\vee$, and is therefore a Poisson map by Theorem
\ref{th:EEM}. 
\end{pf}

\section{Isomonodromic deformations}\label{Section:isomono}

\newcommand {\cl}{^{\operatorname{scl}}}


Let $S_\pm\in\Ug^{\otimes 2}\fml^o$ be the Stokes matrices of the dynamical KZ
equations, and $S_\pm\cl\in G[[\g^*]]_+$ their semiclassical limit.

For any $\alpha\in\sfPhi$, let $Q_\alpha\in S^2\g\subset\IC[\g^*]$
be given by $Q_\alpha=x_\alpha\cdot x_{-\alpha}=Q_{-\alpha}$.\comment
{Work out what the corresponding vector field looks like}

\begin{prop}\hfill\break
\begin{enumerate}
\item As a function of $\mu\in\h\reg^\IR$, $S_\pm\cl$ satisfies the
following PDE
\[d_\h S_\pm\cl=
\frac{1}{2\pi\iota}
\sum_{\alpha\in\sfPhi_+}\frac{d\alpha}{\alpha}\{Q_\alpha,S_\pm\cl\}\]
\item Regard $B\in\g$ as a function of $\mu\in\hreg^\IR$. Then, the
Stokes matrices of (the classical ODE) are locally constant as
$\mu$ varies in $\hreg$ if, and only if $B$ satisfies the nonlinear
PDE
\[d_\h B=
-\frac{1}{2\pi\iota}\sum_{\alpha\in\sfPhi_+}\frac{d\alpha}{\alpha}H_\alpha\]
where $H_\alpha=\{Q_\alpha,-\}$ is the Hamiltonian vector field corresponding
to $Q_\alpha$.
\end{enumerate}
\end{prop}
\begin{pf}
(1)
By Proposition \ref{pr:Spm}, $S_\pm$ satisfy
\[d_\h S_\pm=
\frac{1}{4\pi\iota}\sum_{\alpha\in\sfPhi_+}\frac{d\alpha}{\alpha}
\left[\onetwo{\hbar\Kalpha},S_\pm\right]\]
Note that $\hbar^2\Kalpha\in\UU'$, and that its image in $\UU'/\hbar
\UU'$ is $2Q_\alpha$. As pointed out in \ref{se:filtered A}, $\hbar\ad
(\Kalpha)$ is a derivation of $\UU'$. Since $[\hbar\Kalpha,-]=\hbar^
{-1}[\hbar^2\Kalpha,-]$, $\hbar\ad(\Kalpha)$ induces the derivation
$\{Q_\alpha,-\}$ on $\IC[\g^*]$. The result now follows from the fact
that $\hbar\Kalpha^{(1)}\in\hbar\UU\otimes\UU'$, so that its image
in $\Ug\fmls{\g^*}$ is zero.

(2)
\[\begin{split}
d_\h S_\pm\cl
&=
\frac{1}{2\pi\iota}
\sum_{\alpha\in\sfPhi_+}\frac{d\alpha}{\alpha}\{Q_\alpha,S_\pm\cl\}+
d_{\g^*} S_\pm\cl(d_\h B)\\
&=
d_{\g^*} S_\pm\cl\left(
\frac{1}{2\pi\iota}\sum_{\alpha\in\sfPhi_+}\frac{d\alpha}{\alpha}H_\alpha+d_\h B
\right)\end{split}\]
\end{pf}

This is the time-dependent Hamiltonian description of the isomonodromic deformation given by \cite{JMU, Bo2}. 
Here we give a quantum algebra proof, which enables us to interpret the symplectic nature of the isomonodromic deformation from the perspective of the gauge action of Casimir operators on quantum Stokes matrices.

\end{document}